\input amstex
\magnification=\magstep1 
\baselineskip=13pt
\documentstyle{amsppt}
\vsize=8.7truein \CenteredTagsOnSplits \NoRunningHeads
\def\today{\ifcase\month\or
  January\or February\or March\or April\or May\or June\or
  July\or August\or September\or October\or November\or December\fi
  \space\number\day, \number\year}
 
 \def\EE{\bold{E}\thinspace}
 \def\PP{\bold{P}\thinspace}
 \def\vl{\operatorname{vol}}
 
 \def\spa{\operatorname{span}}
 \def\cov{\bold{cov}\thinspace}
  \def\var{\bold{var}\thinspace}
  \def\dist{\operatorname{dist}}
  \def\rk{\operatorname{rank}}
  
  \def\UU{\Cal U}
  \def\GG{\Cal G}
 \topmatter

\title An asymptotic formula for the number of non-negative integer matrices with prescribed row 
and column sums 
\endtitle
\author Alexander Barvinok and J.A. Hartigan 
\endauthor
\address Department of Mathematics, University of Michigan, Ann Arbor,
MI 48109-1043, USA
\endaddress
\email barvinok$\@$umich.edu \endemail
\address Department of Statistics, Yale University, New Haven, CT 06520-8290 \endaddress
\email john.hartigan$\@$yale.edu \endemail 
\thanks The research of the first author was partially supported by NSF Grant DMS 0856640 and a
United States - Israel BSF grant 2006377. \endthanks
\date April 2010
\enddate
\abstract We count $m \times n$ non-negative integer matrices (contingency tables) 
with prescribed row and column sums (margins). For a wide class of {\it smooth} margins we 
establish a computationally efficient asymptotic formula approximating the number of matrices within 
a relative error which approaches 0 as $m$ and $n$ grow. 
\endabstract
\endtopmatter
\document

\head 1. Introduction and main results \endhead 

Let $R=\left(r_1, \ldots, r_m \right)$ and $C=\left(c_1, \ldots, c_n \right)$ be positive integer vectors 
such that 
$$r_1 + \ldots + r_m = c_1 + \ldots + c_n=N.$$
We are interested in the number $\#(R, C)$ of $m \times n$ non-negative integer matrices $D=\left(d_{ij}\right)$ 
with row sums $R$ and column sums $C$. Such matrices $D$ are often called {\it contingency tables} 
with {\it margins} $(R, C)$. The problem of computing or estimating $\#(R,C)$ efficiently
has attracted considerable attention, see, for example, \cite{B+72}, \cite{Be74}, \cite{GC77}, \cite{DE85},
\cite{DG95}, \cite{D+97}, \cite{Mo02}, \cite{CD03},  \cite{C+05}, \cite{CM07}, \cite{GM07},
\cite{B+08}, \cite{Z+09} and \cite{Ba09}. 

Asymptotic formulas for numbers $\#(R, C)$ as $m$ and $n$ grow are known in sparse cases, where the average entry $N/mn$ of the matrix 
goes to 0, see \cite{B+72}, \cite{Be74}, \cite{GM07} and in the case when 
all row and all column sums are equal, $r_1= \ldots =r_m$ and $c_1=\ldots =c_n$, \cite{CM07}. 
In \cite{Ba09} an asymptotic formula for $\log \#(R,C)$ is established under quite general 
circumstances. 

In this paper, we prove an asymptotic formula for $\#(R, C)$  for a reasonably wide class of 
{\it smooth} margins $(R, C)$. In \cite{BH10} we apply a similar approach to find an 
asymptotic formula for the number of matrices ({\it binary contingency tables}) with row
sums $R$, column sums $C$ and 0-1 entries as well as to find an asymptotic formula 
for the number of graphs with prescribed degrees of vertices.

\subhead (1.1) The typical matrix and smooth margins \endsubhead 
The typical matrix was introduced in \cite{Ba09}  and various versions of smoothness 
for margins were introduced 
in \cite{B+08} and in \cite{Ba08}. The function 
$$g(x) =(x+1) \ln (x+1) - x\ln x \quad \text{for} \quad x \geq 0$$
plays the crucial role. It is easy to see that $g$ is increasing and concave
with $g(0)=0$. For an $m \times n$ non-negative 
matrix $X=\left(x_{jk} \right)$ we define 
$$g(X)=\sum \Sb 1 \leq j \leq m \\ 1 \leq k \leq n \endSb g\left(x_{jk}\right) =
\sum \Sb 1 \leq j \leq m \\ 1 \leq k \leq n \endSb
 \Bigl( \left(x_{jk}+1 \right) \ln \left(x_{jk}+1 \right) -x_{jk} \ln x_{jk} \Bigr).$$

Given margins $(R, C)$, let 
$P(R, C)$ be the polytope of all real non-negative $m \times n$ matrices $X=\left(x_{jk}\right)$ with row 
sums $R$ and column sums $C$, also known as the {\it transportation polytope}.
 We consider the following optimization problem:
$$\text{Find} \quad \max_{X \in P(R, C)} g(X). \tag1.1.1$$
Since $g$ is strictly concave, the maximum is attained at a unique matrix $Z=\left(\zeta_{jk} \right)$, which 
we call the {\it typical} matrix with margins $(R, C)$. One can show that $\zeta_{jk} >0$ for all $j$ and $k$, see \cite{B+08} and \cite{Ba08}. In \cite{Ba08} it is shown that a random contingency 
table, sampled from the uniform distribution on the set of all non-negative integer matrices 
with row sums $R$ and column sums $C$ is, in some rigorously defined sense, likely to be close to the 
typical matrix $Z$. In \cite{BH09} we give the following probabilistic interpretation of $Z$.
Let us consider the family of all probability distributions on the set ${\Bbb Z}_+^{m \times n}$ of all non-negative $m \times n$ integer matrices with the expectations in the affine subspace ${\Cal A}(R, C)$
 of the 
$m \times n$ matrices with row sums $R$ and column sums $C$. In this family there is
a unique distribution of the maximum entropy and $Z$ turns out to be the expectation
of that distribution. The maximum 
entropy distribution is necessarily a distribution on ${\Bbb Z}_+^{m \times n}$ with independent
geometrically distributed coordinates,
which, conditioned on ${\Cal A}(R, C)$, results in the uniform distribution on the set of contingency 
tables with margins $(R, C)$. Function $g(X)$ turns out
to be the entropy of the multivariate geometric distribution on ${\Bbb Z}^{m \times n}_+$ with the 
expectation $X$.

Let us fix a number $0< \delta < 1$. We say that margins $(R, C)$ are $\delta$-{\it smooth} provided
the following conditions (1.1.2)--(1.1.4) are satisfied: 
$$m \ \geq \  \delta n \quad \text{and} \quad n \ \geq \ \delta m, \tag1.1.2$$
so the dimensions of the matrix are of the same order;
$$\delta \tau \ \leq \ \zeta_{jk} \ \leq \ \tau \quad \text{for all} \quad j \quad \text{and} \quad k, \tag1.1.3$$
for some $\tau$ such that
$$\tau \ \geq \ \delta \tag1.1.4$$

We note that $\delta$-smooth margins are also $\delta'$-smooth for any $0<\delta' < \delta$.

Condition (1.1.3) requires that the entries of the typical matrix are of the same order and it plays 
a crucial role in our proofs. Often, one can show that  margins  are smooth by predicting
what the solution to the optimization  problem (1.1.1) will look like. For example, if all row sums $r_j$ are equal, symmetry requires that we have $\zeta_{jk}=c_k/m$ for all $j$ and $k$, so 
the entries of the typical matrix are of the same order provided the column sums $c_k$ are of the same
order. On the other hand, (1.1.3) is violated in some curious cases. For example, if $m=n$ and
$r_1=\ldots =r_{n-1}=c_1=\ldots =c_{n-1}=n$ while $r_n=c_n =3n$, the entry $\zeta_{nn}$ of the 
typical matrix is linear in $n$, namely $\zeta_{nn} > 0.58 n$, while all other entries of $Z$ remain 
bounded by a constant, see \cite{Ba08}. If we change $r_n$ and $c_n$ to $2n$, the entry $\zeta_{nn}$
becomes bounded by a constant as well. One may wonder (this question is inspired by a
conversation with B. McKay) if 
the smoothness condition (1.1.3) is indeed necessary for the number of tables $\#(R, C)$ to be 
expressible by a formula which varies ``smoothly'' as the margins $R$ and $C$ vary, like 
the formula in Theorem 1.3 below.
In particular, can there be a sudden jump in the number of tables with $m=n$, 
$r_1=\ldots =r_{n-1}=c_1= \ldots =c_{n-1}=n$ when $r_n=c_n$ crosses a certain threshold between 
$2n$ and $3n$?

 In \cite{B+08} it is proven that if the ratio of the maximum row sum $r_+=\max_j r_j$ 
to the minimum row sum $r_-=\min_j r_j$ and the ratio of the maximum column sum $c_+=\max_k c_k$
to the minimum column sum $c_-=\min_k c_k$ do not exceed a number  $\beta < (1+\sqrt{5})/2 \approx 1.618$, then (1.1.3) is satisfied with some $\delta=\delta(\beta)>0$. The bound $(1+\sqrt{5})/2$ 
is not optimal, apparently it can be increased to 2, see \cite{Lu08}.
  It looks plausible that if the margins are of the 
same order and sufficiently generic then the entries of the typical table are of the same order 
as well. 
 
The lower bound in (1.1.4) requires that the density $N/mn$ of the margins, that is the average
entry of the matrix, remains bounded away from 0. This is unavoidable as our asymptotic 
formula does not hold for sparse cases where $N/mn \longrightarrow 0$, see \cite{GM07}. 
\bigskip
We proceed to define various objects needed to state our asymptotic formula.

\subhead (1.2) Quadratic form $q$ and related quantities \endsubhead 
Let $Z=\left(\zeta_{jk} \right)$ be the typical matrix defined in Section 1.1. We consider the following quadratic
form $q: {\Bbb R}^{m +n} \longrightarrow {\Bbb R}$:
$$\aligned q\left(s, t \right) ={1 \over 2} \sum \Sb 1 \leq j \leq m \\ 1 \leq k \leq n \endSb 
&\left(\zeta_{jk}^2 +\zeta_{jk} \right) \left(s_j+ t_k \right)^2 \quad \text{where} \\
&s=\left(s_1, \ldots, s_m \right) \quad \text{and} \quad t=\left(t_1, \ldots, t_n \right). \endaligned \tag1.2.1$$
Thus $q$ is a positive semidefinite quadratic form. It is easy to see that the null-space of $q$ is spanned 
by vector 
$$u=\left(\underbrace{1, \ldots, 1}_{\text{$m$ times}}; \underbrace{-1, \ldots, -1}_{\text{$n$ times}} \right).$$
Let $H=u^{\bot}$, $H \subset {\Bbb R}^{m+n}$,
 be the orthogonal complement to $u$. Then the restriction $q|H$ is a positive definite 
quadratic form and hence we can define its determinant $\det q|H$ that is the product of the non-zero eigenvalues 
of $q$. 
Let us define polynomials $f, h: {\Bbb R}^{m+n} \longrightarrow {\Bbb R}$ by 
$$\aligned & f(s, t) ={1 \over 6} \sum \Sb 1 \leq j \leq m \\ 1 \leq k \leq n \endSb \zeta_{jk} \left( \zeta_{jk} +1 \right)
\left(2 \zeta_{jk}+1 \right) \left(s_j +t_k \right)^3 \ \\  &\qquad \text{and} \\
&h(s,t)={1 \over 24} \sum \Sb 1 \leq j \leq m \\ 1 \leq k \leq n \endSb  \zeta_{jk} \left(\zeta_{jk} +1 \right) \left(6 \zeta_{jk}^2 +6\zeta_{jk}+1 \right)
\left(s_j + t_k \right)^4. \endaligned \tag1.2.2$$
We consider the Gaussian probability measure on $H$ with the 
density proportional to $e^{-q}$ and define 
$$\mu = \EE f^2 \quad \text{and} \quad \nu = \EE h.$$
Now we have all the ingredients to state our asymptotic formula for $\#(R, C)$.
\proclaim{(1.3) Theorem} Let us fix $0 < \delta <1$. Let $(R, C)$ be $\delta$-smooth margins,
let the function $g$ and the typical matrix $Z$ be as defined in Section 1.1 and let 
the quadratic form $q$ and values of $\mu$ and $\nu$ be as defined in Section 1.2. Then 
the value of 
$${e^{g(Z)} \sqrt{m+n}  \over (4\pi)^{(m+n-1)/2} \sqrt{\det q|H}} \exp\left\{-{\mu \over 2} + \nu\right\}$$
approximates $\#(R,C)$ within a relative error which approaches 0
as $m, n \longrightarrow +\infty$. More precisely, for any $0< \epsilon \leq 1/2$ the above expression
approximates $\#(R,C)$ within relative error $\epsilon$ provided 
$$m+n \ \geq \ \left({1\over \epsilon}\right)^{\gamma(\delta)},$$
for some $\gamma(\delta)>0$.
\endproclaim

In \cite{CM07} Canfield and McKay obtain an asymptotic formula for $\#(R, C)$
in the particular case of all row sums being equal and all column sums being equal. 
One can show that our formula indeed becomes the asymptotic formula of \cite{CM07}
when $r_1=\ldots =r_m$ and $c_1=\ldots =c_n$.
 In \cite{Ba09} it is proven that  the value $g(Z)$
provides an asymptotic approximation to $\ln \#(R, C)$ for a rather wide class of margins
(essentially, we need only the density $N/mn$ to be bounded away from 0 but do not need a
subtler condition (1.1.3) of smoothness). The first part 
$${e^{g(Z)} \sqrt{m+n}  \over (4\pi)^{(m+n-1)/2} \sqrt{\det q|H}} \tag1.3.1$$
of the formula is called the ``Gaussian approximation'' in \cite{BH09}. It has the following intuitive 
explanation. Let us consider a random matrix $X$ with the multivariate geometric distribution on 
the set ${\Bbb Z}^{m \times n}_+$ of all non-negative integer matrices such that $\EE X=Z$,
where $Z$ is the typical matrix with margins $(R, C)$. It follows from the results of \cite{BH09} that
the distribution of $X$ conditioned on the affine subspace ${\Cal A}={\Cal A}(R,C)$ of matrices with row sums 
$R$ and column sums $C$ is uniform with the probability mass function of $e^{-g(Z)}$ for 
every non-negative integer matrix in ${\Cal A}$.
Therefore, 
$$\#(R, C)=e^{g(Z)} \PP\bigl\{X \in {\Cal A} \bigr\}.$$
Let $Y \in {\Bbb R}^{m+n}$ be a random vector obtained by computing $m$ row sums and $n$ 
column sums of $X$. Then $\EE Y=(R,C)$ and
$$\PP\bigl\{X \in {\Cal A} \bigr\} =\PP\bigl\{Y=(R, C)\bigr\}.$$
We obtain (1.3.1) if we assume in the spirit of the Local Central Limit Theorem that the distribution 
of $Y$ in the vicinity of $\EE Y$ is close to the $(m+n-1)$-dimensional Gaussian distribution
(we lose one dimension since the row and column sums of a matrix are bound by one linear 
relation: the sum of all row sums is equal to the sum of all column sums). This assumption is 
not implausible since the coordinates of $Y$ are obtained by summing up of a number of independent 
entries of $X$. 

The correction factor
$$\exp\left\{-{\mu \over 2} +\nu \right\} \tag1.3.2$$
is, essentially, the {\it Edgeworth correction} in the Central Limit Theorem. In the course of the proof
of Theorem 1.3 we establish a two-sided bound
$$\gamma_1(\delta) \ \leq \ \exp\left\{-{\mu \over 2} +\nu \right\} \ \leq \ \gamma_2(\delta)$$
for some constants $\gamma_1(\delta), \gamma_2(\delta)>0$ as long as the margins $(R, C)$
remain $\delta$-smooth.

De Loera \cite{D09a}, \cite{D09b} ran a range of numerical experiments which seem to demonstrate that 
already the Gaussian approximation (1.3.1) works reasonably well for contingency tables.
For example, for $R=\left(220, 215, 93, 64\right)$ and $C=\left(108, 286, 71, 127 \right)$ 
formula (1.3.1) approximates $\#(R, C)$ within a relative error of about $6\%$, for 
$R=C=\left(300, 300, 300, 300 \right)$ the error is about $12\%$ while for 
$R=\left(65205, 189726, 233525, 170004\right)$ and $C=\left(137007, 87762, 274082, 159609\right)$
the error is about $1.2\%$. 

\subhead (1.4) Computations and a change of the hyperplane \endsubhead 
Optimization problem (1.1.1) is convex and can be solved, for example, by interior point methods,
see \cite{NN94}. That is, 
for any $\epsilon>0$ the entries $\zeta_{jk}$ of the typical matrix $Z$ can be computed within relative
error $\epsilon$ in time polynomial in $\ln (1/\epsilon)$ and $m+n$.

Given $Z$, quantities $\det q|H$, $\mu$ and $\nu$ can be computed by linear algebra algorithms
in $O\left( m^2 n^2 \right)$ time, 
since to compute the expectation of a polynomial with respect to the Gaussian measure one only needs 
to know the covariances of the variables, see Section 4.2. It may be convenient to replace 
the hyperplane $H \subset {\Bbb R}^{m+n}$ orthogonal to the null-space of $q$ by a coordinate 
hyperplane $L \subset {\Bbb R}^{m+n}$ defined by any of the equations $s_j=0$ or $t_k=0$.
Indeed, if $L \subset {\Bbb R}^{m+n}$ is any hyperplane not containing the null-space of $q$, then 
the restriction $q|L$ is strictly positive definite and one can consider the Gaussian probability 
measure in $L$ with the density proportional to $e^{-q}$. We prove in Lemma 3.1 below that 
the expectation of any polynomial in $s_j +t_k$ does not depend on the choice of $L$ and hence 
$\mu$ and $\nu$ can be defined as in Section 1.2 with $H$ replaced by $L$. We describe the dependence 
of $\det q|L$ on $L$ in Lemma 3.5. In particular, it follows that if $L$ is 
a coordinate hyperplane then $\det q|H =(m+n) \det q|L$. 

If we choose $L$ to be defined 
by the equation $t_n=0$ then we have an explicit formula for the matrix $Q$ of $q|L$ as follows:
$$q(x) = {1 \over 2} \langle x, Qx \rangle \quad \text{for} \quad 
x=\left(s_1, \ldots, s_m; t_1, \ldots, t_{n-1} \right),$$
where $\langle \cdot, \cdot \rangle$ is the standard scalar product and 
$Q=\left(q_{il}\right)$ is the $(m+n-1) \times (m+n-1)$ symmetric matrix, where 
$$\split q_{j(k+m)}=q_{(k+m)j}=&\zeta_{jk}^2 +\zeta_{jk} \quad \text{for} \quad j=1, \ldots, m 
\quad  \text{and} \quad k=1, \ldots, n-1, \\
q_{jj}=&r_j + \sum_{k=1}^n \zeta_{jk}^2  \quad \text{for} \quad j=1, \ldots, m, \\ 
q_{(k+m)(k+m)}=&c_k + \sum_{j=1}^n \zeta_{jk}^2 \quad \text{for} \quad k=1, \ldots, n-1,\endsplit$$
and all other entries $q_{il}$ are zeros.
Then $\det q|L = 2^{1-m-n} \det Q$ and $Q^{-1}$ is the covariance matrix of $s_1, \ldots, s_m; 
t_1, \ldots, t_{n-1}$.

\head 2. An integral representation for the number of contingency tables and the plan 
of the proof of Theorem 1.3\endhead 

In \cite{BH09} we prove the following general result.
\proclaim{(2.1) Theorem} Let $P \subset {\Bbb R}^p$ be a polyhedron defined by the system of linear
equations $Ax=b$, where $A$ is a $d \times p$ integer matrix with columns 
$a_1, \ldots, a_p \in {\Bbb Z}^d$ and $b \in {\Bbb Z}^d$ is an integer vector, and inequalities 
$x \geq 0$ (the inequalities are understood as coordinate-wise). Suppose that $P$ is bounded 
and has a non-empty interior, that is, contains a point $x=\left(\xi_1, \ldots, \xi_p \right)$ 
such that $\xi_j >0$ for $j=1, \ldots, p$. Then the function
$$g(x)=\sum_{j=1}^p \Bigl( \left(\xi_j +1 \right) \ln \left( \xi_j +1 \right) - \xi_j \ln \xi_j \Bigr)$$
attains its maximum on $P$ at a unique point $z=\left(\zeta_1, \ldots, \zeta_p\right)$ such that 
$\zeta_j >0$ for $j=1, \ldots, p$.

Let $\Pi \subset {\Bbb R}^d$ be the parallelepiped consisting of the points 
$t=\left(\tau_1, \ldots, \tau_d\right)$ such that 
$$ -\pi \ \leq \ \tau_k \ \leq \ \pi \quad \text{for} \quad k=1, \ldots, d.$$
Then the number $|P \cap {\Bbb Z}^p|$ of integer points in $P$ can be written as 
$$|P \cap {\Bbb Z}^p| ={e^{g(z)} \over (2 \pi)^d} \int_{\Pi} e^{-i \langle t, b \rangle} 
\prod_{j=1}^p {1 \over 1+\zeta_j - \zeta_j e^{i \langle a_j, t \rangle}} \ d t,$$
where $\langle \cdot, \cdot \rangle$ is the standard scalar product in ${\Bbb R}^d$ and $i=\sqrt{-1}$.
\endproclaim 
{\hfill \hfill \hfill} \qed

The idea of the proof is as follows. Let $X =\left(x_1, \ldots, x_p \right)$ be a random 
vector of independent geometric random variables $x_j$ such that $\EE x_j =\zeta_j$.
Hence values of $X$ are non-negative integer vectors and we show in \cite{BH09} that the 
probability mass function of $X$ is constant on the set $P \cap {\Bbb Z}^p$ and equals $e^{-g(z)}$ 
for every integer point in $P$. Letting $Y=AX$, we obtain
$$|P\cap {\Bbb Z}^p| =e^{g(z)} \PP\{X \in P\} =e^{g(z)} \PP\{ Y =b\}$$ 
and the probability in question is written as the integral of the characteristic function of $Y$.

Since 
$$\sum_{j=1}^p \zeta_j a_j =b,$$
in a neighborhood of the origin $t=0$ the integrand can be written
as 
$$\aligned 
e^{-i \langle t, b \rangle} \prod_{j=1}^p & {1 \over 1+\zeta_j - \zeta_j e^{i \langle a_j, t \rangle}} 
\\ =\exp\Biggl\{
&-{1 \over 2} \sum_{j=1}^p  \left(\zeta_j^2 + \zeta_j \right) \langle a_j, t \rangle^2 \\
&-{i \over 6}  \sum_{j=1}^p \zeta_j \left(\zeta_j +1 \right) \left(2\zeta_j+1 \right) 
\langle a_j, t \rangle^3 \\ &+{1 \over 24} \sum_{j=1}^p \zeta_j \left(\zeta_j+1\right) \left(6 \zeta_j^2 +
6 \zeta_j+1 
\right) \langle a_j, t \rangle^4 \\&+O\left( \sum_{j=1}^p \left(\zeta_j +1 \right)^5 \langle a_j, t \rangle^5
\right) \Biggr\}. \endaligned \tag2.1.1$$
Note that the linear term is absent in the expansion.

We obtain the following corollary. 
\proclaim{(2.2) Corollary} Let $R=\left(r_1, \ldots, r_m\right)$ and $C=\left(c_1, \ldots, c_n \right)$
be margins and let $Z=\left(\zeta_{jk}\right)$ be the typical matrix defined in Section 1.1.
Let
$$F(s,t) =\exp\left\{-i \sum_{j=1}^m r_j s_j -i\sum_{k=1}^n c_k t_k \right\}
\prod \Sb 1 \leq j \leq m \\ 1 \leq k \leq n \endSb {1 \over 1+\zeta_{jk} - \zeta_{jk} 
e^{i\left(s_j +t_k\right)}}.$$
Let $\Pi \subset {\Bbb R}^{m+n}$ be the parallelepiped consisting of the points 
$\left(s_1, \ldots, s_m; t_1, \ldots, t_n \right)$ such that 
$$-\pi \ \leq \ s_j, t_k \ \leq \ \pi \quad \text{for all} \quad j, k.$$
Let us identify ${\Bbb R}^{m+n-1}$ with the hyperplane $t_n=0$ in ${\Bbb R}^{m+n-1}$ 
and let $\Pi_0 \subset \Pi$ be the facet of $\Pi$ defined by the equation $t_n=0$.
Then
$$\#(R,C)={e^{g(Z)} \over (2 \pi)^{m+n-1}} \int_{\Pi_0} F(s,t) \ ds dt,$$
where $ds dt$ is the Lebesgue measure in $\Pi_0$.
\endproclaim
\demo{Proof} The number $\#(R,C)$ of non-negative integer $m \times n$ matrices with row sums $R$ and column sums $C$ is the number of integer points in the transportation polytope $P(R,C)$.
 We can define $P(R,C)$ by prescribing all 
row sums $r_1, \ldots, r_m$ and all but one column sums $c_1, \ldots, c_{n-1}$ of a non-negative 
$m \times n$ matrix.
Applying Theorem 2.1, we get the desired integral representation.
{\hfill \hfill \hfill} \qed
\enddemo

From (2.1.1) we get the following expansion in the neighborhood of 
$s_1 =\ldots = s_m=t_1 =\ldots =t_n=0$:
$$\aligned F(s, t) =\exp&\Biggl\{ -q(s, t) -i f(s, t) +h(s, t) \\
&+ O\left(\sum_{j,k} \left(1 +\zeta_{jk}\right)^5 
\left(s_j +t_k \right)^5 \right) \Biggr\}, \endaligned  \tag2.2.1 $$
where $q$, $f$, and $h$ are defined by (1.2.1)--(1.2.2).

\subhead (2.3) The plan of the proof of Theorem 1.3 \endsubhead First, we argue that it suffices 
to prove Theorem 1.3 under one additional assumption, namely, that the parameter $\tau$ in 
(1.1.3) is bounded by a polynomial in $m+n$:
$$\tau \ \leq \ (m+n)^{1/\delta} \quad \text{for some} \quad \delta >0 $$
(for example, one can choose $\delta=1/10$). Indeed, it follows by results of \cite{D+97} (see Lemma 3
there) that for $\tau \geq (mn)^2$ the (properly normalized) volume $\vl P(R, C)$ of the transportation 
polytope approximates the number of tables $\#(R, C)$ within 
a relative error of $O\left((m+n)^{-1}\right)$. Since $\dim P(R, C)=(m-1)(n-1)$ and 
$$\vl P(\alpha R, \alpha C) = \alpha^{(m-1)(n-1)} \vl P(R, C) \quad \text{for} \quad \alpha >0,$$
to handle larger $\tau$ it suffices to show that the formula of Theorem 1.3 scales the right way
if the margins $(R, C)$ get scaled $(R, C) \longmapsto (\alpha R, \alpha C)$ (and appropriately 
rounded, if the obtained margins are not integer). If $\tau$ is large enough then scaling results 
in an approximate scaling $e^{g(Z)} \longmapsto \alpha^{mn} e^{g(Z)}$, $q \longmapsto \alpha^2 q$, 
$f \longmapsto \alpha^3 f$ and $h \longmapsto \alpha^4 h$ and hence the value produced by 
the formula of Theorem 1.3 gets multiplied by roughly $\alpha^{(m-1)(n-1)}$, as desired.
We provide necessary details in Section 8.

To handle the case of $\tau$ bounded by a polynomial in $m+n$, we use the integral representation of Corollary 2.2. 

Let us define a neighborhood $\UU \subset \Pi_0$ of the origin by 
$$\UU=\left\{ \left(s_1, \ldots, s_m; t_1, \ldots, t_{n-1} \right): \quad \left|s_j\right|, \left|t_k \right| 
\ \leq \ {\ln (m+n) \over \tau \sqrt{m+n}} \quad \text{for all} \quad j, k \right\}.$$
We show that the integral of $F(s, t)$ over $\Pi_0 \setminus \UU$ is asymptotically negligible. 
Namely, in Section 7 we prove that the integral 
$$\int_{\Pi_0 \setminus \UU} |F(s, t)| \ ds dt $$
is asymptotically negligible compared to the integral 
$$\int_{\UU} |F(s, t)| \ ds dt. \tag2.3.1$$
In Section 6, we evaluate the integral 
$$\int_{\UU} F(s, t) \ ds dt \tag2.3.2$$ 
and show that it produces the asymptotic formula of Theorem 1.3. In particular, we show that 
(2.3.1) and (2.3.2) are of the same order, that is,
$$ \int_{\UU} |F(s, t)| \ ds dt \ \leq \ \gamma(\delta) \left| \int_{\UU} F(s, t) \ ds dt \right|$$
for some constant $\gamma(\delta) \geq 1$. Hence the integral of $F(s, t)$ outside of $\UU$ is 
indeed asymptotically irrelevant.

From (2.2.1), we deduce that 
$$F(s, t) \approx \exp\left\{-q(s, t) -i f(s, t) +h(s, t) \right\} \quad \text{for} \quad (s, t) \in \UU,$$
where $q$ is defined by (1.2.1) and $f$ and $h$ are defined by (1.2.2), so that the 
contribution of the terms 
of order 5 and higher in (2.2.1) is asymptotically negligible in the integral (2.3.2).
The integral of $e^{-q}$ over $\UU$ produces the Gaussian term (1.3.1)
 However, both the cubic term $f(s, t)$ and the fourth-order term $h(s, t)$ contribute substantially to the 
integral, correcting the Gaussian term (1.3.1) by a constant factor.

Let us consider the Gaussian probability measure in the coordinate hyperplane $t_n=0$, 
which we identify with ${\Bbb R}^{m+n-1}$,
with the density proportional to $e^{-q}$. 
 In Section 5, we show that with respect to that measure, $h(s, t)$ remains, essentially, constant in the neighborhood $\UU$:
$$h(s, t) \approx \EE h = \nu \quad \text{almost everywhere in} \quad \UU.$$
This allows us to conclude that asymptotically
$$\int_{\UU} \exp\bigl\{-q(s, t)-if(s, t) +h(s, t) \bigr\} \ ds dt \approx 
e^{\nu} \int_{\UU} \exp\left\{-q(s, t) +i f(s, t)\right\} \ ds dt.$$

In Section 4, we show that $f(s, t)$ behaves, essentially, as a Gaussian random variable with 
respect to the probability measure in ${\Bbb R}^{m+n-1}$ with the density proportional to $e^{-q}$, so
$$\split \int_{\UU} \exp \bigl\{-q(s, t)-if(s, t)\bigr\} \ ds dt \approx &\int_{{\Bbb R}^{m+n-1}} \exp\bigl\{-q(s, t) -
if(s, t) \bigr\} \ ds dt \\ \approx &\exp\left\{ -{1 \over 2} \EE f^2 \right\} \int_{{\Bbb R}^{m+n-1}} e^{-q(s, t)} 
\ ds dt, \endsplit$$
which concludes the computation of (2.3.2).

The results of Sections 4 and 5 are based on the analysis in Section 3. In Section 3, we consider 
coordinate functions $s_j$ and $t_k$ as random variables with respect to the Gaussian 
probability measure on a hyperplane $L \subset {\Bbb R}^{m+n}$ not containing the null-space of $q$
with the density proportional to $e^{-q}$.
We show that $s_{j_1} +t_{k_1}$ and $s_{j_2} +t_{k_2}$ are weakly correlated
provided $j_1 \ne j_2$ and $k_1 \ne k_2$, that is,
$$\split \EE  \left| \left(s_{j_1} +t_{k_1} \right) \left(s_{j_2} +t_{k_2}\right) \right| =&O \left({1 \over mn}\right)
\quad \text{provided} \quad j_1 \ne j_2 \quad \text{and} \quad k_1 \ne k_2 \quad \text{and} \\
\EE \left| \left(s_j+t_k \right)^2 \right| =&O\left({1 \over m+n} \right) \quad \text{for all} \quad j, k.
\endsplit$$
\subhead (2.4) Notation \endsubhead In what follows, we denote by $\gamma$, sometimes with 
an index or a list of parameters, a positive constant depending on the parameters. The actual 
value of $\gamma$ may change from line to line. The most common appearance will be 
$\gamma(\delta)$, a positive constant depending only on the $\delta$-smoothness constant 
$\delta$.

As usual, for two functions $f$ and $g$, where $g$ is non-negative, we say that $f=O(g)$ if $|f| \leq \gamma g$ for 
some constant $\gamma >0$ and that $f =\Omega(g)$ if $f \geq \gamma g$ for some constant 
$\gamma>0$.

\head 3. Correlations \endhead

Recall (see Section 1.2) that the quadratic form  $q: {\Bbb R}^{m+n} \longrightarrow {\Bbb R}$ is
defined by 
$$q(s, t)={1 \over 2} \sum \Sb 1 \leq j \leq m \\ 1 \leq k \leq n \endSb \left(\zeta_{jk}^2 +\zeta_{jk}\right) 
\left(s_j +t_k \right)^2 \quad \text{for} \quad (s, t)=\left(s_1, \ldots, s_m; t_1, \ldots t_n \right).$$
Let 
$$u=\left(\underbrace{1, \ldots, 1}_{\text{$m$ times}};\ \underbrace{-1, \ldots, -1}_{\text{$n$ times}} \right).$$
Let $L \subset {\Bbb R}^{m+n}$ be a hyperplane which does not contain $u$. Then the 
restriction $q|L$ of $q$ onto $L$ is a positive definite quadratic form and we can consider the 
Gaussian probability measure on $L$ with the density proportional to $e^{-q}$. We consider $s_j$ and $t_k$ as random variables on $L$ and estimate their covariances.

\proclaim{(3.1) Lemma} For any $1 \leq j_1, j_2 \leq m$ and any $1 \leq k_1, k_2 \leq n$ the 
covariance 
$$\EE \left(s_{j_1} + t_{k_1}\right) \left(s_{j_2} + t_{k_2} \right)$$ 
is independent on the choice the hyperplane $L$, as long as $L$ does not contain $u$.
\endproclaim
\demo{Proof} Let $L_1, L_2 \subset {\Bbb R}^n$ be two hyperplanes not containing $u$. Then 
we can define the projection $pr: L_1 \longrightarrow L_2$ along the span of $u$, so that 
$pr(x)$ for $x \in L_1$ is the unique $y \in L_2$ such that $y-x$ is a multiple of $u$. 
We note that $q(x)=q(x + t u)$ for all $x \in {\Bbb R}^{m+n}$ and all $t \in {\Bbb R}$.
Therefore, the push-forward of the Gaussian probability measure on $L_1$ with the density proportional
to $e^{-q}$ is the probability measure on $L_2$ with the density proportional to $e^{-q}$. 
Moreover, the value of $s_j +t_k$ does not change under the projection and hence the result
follows.
{\hfill \hfill \hfill} \qed
\enddemo

The main result of this section is the following theorem.
\proclaim{(3.2) Theorem} Let us fix a number $\delta >0$ and suppose that 
$$\tau \delta \ \leq \ \zeta_{jk} \ \leq \ \tau \quad \text{for all} \quad j, k$$
and some $\tau >0$. 
Suppose that $\delta m \leq n$ and $\delta n \leq m$.

Let us define 
$$\split &a_j=\sum_{k=1}^n \left(\zeta_{jk}^2 + \zeta_{jk} \right) \quad \text{for} \quad j=1, \ldots, m 
\quad \text{and} \\
&b_k =\sum_{j=1}^m \left(\zeta_{jk}^2 + \zeta_{jk} \right) \quad \text{for} \quad k=1, \ldots, n.
\endsplit$$

Let 
$$\Delta = {12 \over \delta^{15/2}  \left(\tau^2 + \tau\right) mn}.$$

Let $L \subset {\Bbb R}^{m+n}$ be a hyperplane not containing the null-space of $q$.
Let us consider the Gaussian probability measure on $L$ with the density proportional 
to $e^{-q}$.

Then 
$$\split &\Big| \EE  \left(s_{j_1} + t_{k_1} \right) \left(s_{j_2} + t_{k_2} \right) \Big| 
 \ \leq \ \Delta  \quad  \text{provided} \quad 
 j_1 \ne j_2 \quad \text{and} \quad k_1 \ne k_2, \\
   &\left| \EE \left(s_j + t_{k_1} \right) \left(s_j + t_{k_2} \right)- {1 \over a_j} \right| 
 \ \leq \ \Delta  \quad  \text{provided}  \quad k_1 \ne k_2, \\
 &\left| \EE  \left(s_{j_1} + t_k \right) \left(s_{j_2} + t_k \right) -{1 \over b_k} \right| 
 \ \leq \ \Delta  \quad  \text{provided} \quad 
 j_1 \ne j_2 \quad \text{and} \\
 &\left| \EE \left(s_j + t_k \right)^2  -{1 \over a_j} -{1 \over b_k} \right| 
 \ \leq \ \Delta  \quad  \text{for all} \quad 
j  \quad \text{and} \quad k.
 \endsplit$$
\endproclaim

The gist of Theorem 3.2 is that for a fixed $\delta>0$, while generally the covariance of $s_{j_1}+t_{k_1}$ and $s_{j_2}+t_{k_2}$ is $O\left( \tau^{-2} (m+n)^{-1}\right)$, it is only 
$O\left(\tau^{-2} (m+n)^{-2} \right)$ when $j_1 \ne j_2$ and $k_1 \ne k_2$.

 In what follows, we will often deal with the following 
situation. Let $V$ be Euclidean space, let $\phi: V \longrightarrow {\Bbb R}$ be a positive 
semidefinite quadratic form and let $L \subset V$ be a subspace such that the restriction 
of $\phi$ onto $L$ is strictly positive definite. We consider the Gaussian probability 
measure on $L$ with the density proportional to $e^{-\phi}$. For a polynomial
(random variable) $f: L \longrightarrow {\Bbb R}$ we denote by $\EE (f;\ \phi | L)$
the expectation of $f$ with respect to that Gaussian measure. Instead of $\EE(f;\ \phi|V)$ we write simply 
$\EE(f;\ \phi)$. 

We will use the following standard facts. 
Suppose that there is a direct sum decomposition $V =L_1 + L_2 + \ldots + L_k$ where 
$L_i$ are pairwise orthogonal, 
such that 
$$\phi\left(x_1 + \ldots + x_k \right) =\sum_{i=1}^k \phi_i \left(x_i \right) \quad \text{for all} \quad 
x_i \in L_i.$$
In other words, the components  $x_i \in L_i$ of a random point $x=x_1 + \ldots + x_k$, $x \in V$, are independent.
Then for any two linear functions $\ell_1, \ell_2:\  V \longrightarrow {\Bbb R}$ we have 
$$\EE \left(\ell_1 \ell_2; \ \phi \right) =\sum_{i=1}^k \EE \left(\ell_1 \ell_2; \ \phi|L_i \right).$$
Indeed, since 
$$\ell_{1,2}\left(x_1 + \ldots + x_k \right)=\sum_{i=1}^k \ell_{1,2}\left(x_i \right),$$
we obtain
$$\EE \left(\ell_1 \ell_2; \ \phi \right)=\sum_{i_1=1}^k \sum_{i_2=1}^k 
\EE \left( \ell_1\left(x_{i_1}\right) \ell_2\left(x_{i_2}\right); \ \phi \right).$$
If $i_1 \ne i_2$, we have 
$$\EE \left(\ell_1\left(x_{i_1}\right) \ell_2\left(x_{i_2}\right); \ \phi \right)=
\EE \left(\ell_1; \ \phi| L_{i_1} \right) \EE \left(\ell_2; \ \phi|L_{i_2} \right)=0$$
while for $i_1=i_2=i$ we have 
$$\EE \left(\ell_1\left(x_i \right) \ell_2\left(x_i \right); \ \phi \right) =
\EE \left( \ell_1 \ell_2; \ \phi|L_i \right).$$

We deduce Theorem 3.2 from the following statement.

\proclaim{(3.3) Proposition} Let $m$ and $n$ be positive integers such that 
$$\delta m \ \leq \ n \quad \text{and} \quad \delta n \ \leq  \ m \quad \text{for some} \quad 0 < \delta < 1.$$
Let  $\xi_{jk}$, $j=1, \ldots, m$ and $k=1, \ldots, n$, 
be real numbers such that 
$$\alpha \ \leq \ \xi_{jk} \ \leq \ \beta \quad \text{for all} \quad j, k$$
and some $\beta > \alpha > 0$.
 Let
$$\split &a_j =\sum_{k=1}^n \xi_{jk} \quad \text{for} \quad j=1, \ldots, m \quad \text{and} \\
&b_k=\sum_{j=1}^m \xi_{jk} \quad \text{for} \quad k=1, \ldots, n. \endsplit$$
Let us define a quadratic form 
$\psi: {\Bbb R}^{m+n} \longrightarrow {\Bbb R}$ by 
$$\split \psi(s, t) = {1  \over 2} \sum \Sb 1 \leq j \leq m \\ 1 \leq k \leq n \endSb &\xi_{jk} \left({s_j \over \sqrt{a_j}} +
{t_k \over \sqrt{b_k}} \right)^2 \\ \text{for} \quad &(s, t)=\left(s_1, \ldots, s_m; t_1, \ldots, t_n\right).
\endsplit$$
Let $L \subset {\Bbb R}^{m+n}$ be the hyperplane consisting of the points 
$\left(s_1, \ldots, s_m; \ t_1, \ldots, t_n \right)$ such that 
$$ \sum_{j=1}^m s_j  \sqrt{a_j}\ =\ \sum_{k=1}^n t_k  \sqrt{b_k}.$$
Then the restriction $\psi|L$ of $\psi$ onto $L$ is strictly positive definite and
for
$$\Delta=3 \left({\beta \over \alpha} \right)^{7/2} {1 \over \sqrt{\delta mn}}$$ we have 
$$ \split &\left|\EE \left(s_j^2;\ \psi|L \right) -1 \right|, \  
\left| \EE \left(t_k^2;\ \psi|L \right) -1 \right| \ \leq \ \Delta \quad \text{for all} \quad j, k, \\
&\left| \EE \left(s_{j_1} s_{j_2};\ \psi|L \right) \right|, \  \left| \EE \left(t_{k_1} t_{k_2};\ \psi|L \right) \right|  \ \leq \  \Delta  \quad \text{for all} \quad j_1 \ne j_2 \quad \text{and} \quad k_1 \ne k_2, \\
& \left| \EE \left(s_j t_k;\ \psi|L \right) \right| \ \leq \ \Delta \quad \text{for all} \quad j, k. \endsplit$$
\endproclaim
\demo{Proof} Clearly, the null-space of $\psi$ is one-dimensional and spanned by vector 
$$w=\left(\sqrt{a_1}, \ldots, \sqrt{a_m}; \ -\sqrt{b_1}, \ldots, -\sqrt{b_n} \right).$$
We have 
$L=w^{\bot}$ and hence the restriction of $\psi$ onto $L$ is positive definite.

Next, we observe that 
$$v=\left( \sqrt{a_1}, \ldots, \sqrt{a_m}; \ \sqrt{b_1}, \ldots, \sqrt{b_n} \right)$$
is an eigenvector of $\psi$ with eigenvalue 1. Indeed, the gradient of $\psi(x)$ at $x=v$ is 
equal to $2v$:
$$\split &{\partial \over \partial s_j} \psi \Big|_{s_j = \sqrt{a_j}, t_k=\sqrt{b_k}} \ = \ {2 \over \sqrt{a_j}} 
\sum_{k=1}^n \xi_{jk} \ = \ 2 \sqrt{a_j} \quad \text{and} \\
 &{\partial \over \partial t_k} \psi \Big|_{s_j = \sqrt{a_j}, t_k=\sqrt{b_k}}\ = \ {2 \over \sqrt{b_k}} 
\sum_{j=1}^m \xi_{jk} \ =\ 2 \sqrt{b_k}. \endsplit$$ 
We write 
$$\psi(s, t) ={1 \over 2} \sum_{j=1}^m s_j^2 \ + \ {1 \over 2} \sum_{k=1}^n t_k^2\
 +\  \sum \Sb j=1, \ldots, m \\ k=1, \ldots, n \endSb
{\xi_{jk} \over \sqrt{a_j} \sqrt{b_k}} s_j t_k. \tag3.3.1$$
Let 
$$c=a_1 + \ldots + a_m =b_1 + \ldots + b_n$$
and let us consider another quadratic form $\phi: {\Bbb R}^{m+n} \longrightarrow {\Bbb R}$ defined by 
$$\phi(s, t) =  \sum \Sb 1 \leq j \leq m \\ 1 \leq k \leq n \endSb {\sqrt{a_j b_k} \over c} s_j t_k 
\qquad \text{for} \quad (s, t)=\left(s_1, \ldots, s_m; t_1, \ldots, t_n \right).  \tag3.3.2$$
Clearly, $\phi(s, t)$ is a form of rank 2. Its non-zero eigenvalues are $-1/2$ with the 
eigenspace spanned by $w$ and $1/2$ with the eigenspace spanned by $v$.

Let us define a subspace $L_0 \subset {\Bbb R}^{m+n}$ of codimension 2 by 
$$L_0 =(v, w)^{\bot}.$$
In other words, $L_0$ consists of the points $\left(s_1, \ldots, s_m; t_1, \ldots, t_n \right)$
such that
$$\sum_{j=1}^m s_j \sqrt{a_j} \ =\ \sum_{k=1}^n t_k \sqrt{b_k}=0.$$
In particular, 
$$\phi(s, t) =0 \quad \text{for all}  \quad (s, t) \in L_0.$$

Let us define a quadratic form 
$$\tilde{\psi} = \psi - \epsilon^2 \phi  \quad \text{for} \quad \epsilon ={\alpha \over \beta}. \tag3.3.3$$
We note that $\tilde{\psi}$ is strictly positive definite. Indeed, $w$ and $v$ are eigenvectors
of $\tilde{\psi}$ with the eigenvalues $\epsilon^2/2 >0$ and $1-\epsilon^2/2 >0$ respectively 
and $\tilde{\psi}$ coincides with $\psi$ on the subspace $L_0 =(v, w)^{\bot}$, where $\psi$ 
is positive definite. Our immediate goal is to bound the covariances
$$
\EE\left(s_{j_1} s_{j_2};\ \tilde{\psi} \right), \ \EE\left(t_{k_1} t_{k_2};\ \tilde{\psi} \right) \quad 
\text{and} \quad \EE \left(s_j t_k;\ \tilde{\psi} \right).$$
We can write 
$$\tilde{\psi}(x) ={1 \over 2} \langle x,\ (I +P) x \rangle \quad \text{for} \quad x=(s, t),$$
where $I$ is the $(m+n)\times (m+n)$ identity matrix, 
$P=\left(p_{il} \right)$ is a symmetric \break $(m+n) \times (m+n)$ matrix with zero 
diagonal and $\langle \cdot, \cdot \rangle$
is the standard scalar product in ${\Bbb R}^{m+n}$.
Since 
$$\aligned  \alpha n \ \leq \ &a_j \ \leq \ \beta n \quad \text{for} \quad j=1, \ldots, m \\
              \alpha m \ \leq \ &b_k \ \leq \ \beta m \quad \text{for} \quad k=1, \ldots, n \quad \text{and} \\
              &c \ \geq \ \alpha mn, \endaligned \tag3.3.4$$
by (3.3.1) -- (3.3.3), for the entries $p_{il}$ of $P$ we have
$$0 \ \leq \ p_{il} \ \leq \ {\beta \over \alpha \sqrt{mn}}={1 \over \epsilon \sqrt{mn}}
 \quad \text{for all} \quad i,l. \tag3.3.5$$
Furthermore, $v$ is the Perron-Frobenius eigenvector of $P$ 
with the corresponding eigenvalue $1-\epsilon^2$.

Let us bound the entries of a positive integer power $P^d=\left(p_{il}^{(d)} \right)$ of $P$.
Let 
$$\kappa ={\beta \over \alpha^{3/2} \delta^{1/4} (mn)^{3/4}} \quad  \text{and let } \quad y=\kappa v, \quad
y=\left(\eta_1, \ldots, \eta_{m+n} \right).$$
From (3.3.4) we conclude that 
$$ a_j, \ b_k \ \geq \ \alpha \sqrt{\delta mn} \quad \text{for all} \quad j, k$$ and hence
by (3.3.5) 
$$ p_{il} \ \leq \ \eta_i \quad \text{for all} \quad i, l. \tag3.3.6$$
Similarly, from (3.3.4), we conclude 
$$a_j, \ b_k \ \leq \ \beta \sqrt{mn/\delta} \quad \text{for all} \quad j, k$$
and hence 
$$\eta_i \ \leq \  {\beta^{3/2} \over \alpha^{3/2} \sqrt{\delta mn}} ={1 \over \epsilon^{3/2} \sqrt{\delta mn}} \quad \text{for all} \quad i. \tag3.3.7$$
Besides, $y$ is an eigenvector of $P^d$ with the eigenvalue $(1-\epsilon^2)^d$.
Therefore, for $d \geq 0$ we have 
$$\split p_{il}^{(d+1)} =&\sum_{j=1}^{m+n} p_{ij}^{(d)} p_{jl} \\ \leq &\sum_{j=1}^{m+n} p_{ij}^{(d)} \eta_j
=(1-\epsilon^2)^d \eta_i  \\ \leq &(1-\epsilon^2)^d {1 \over \epsilon^{3/2} \sqrt{\delta mn}}.
\endsplit$$
Consequently, the series 
$$(I+P)^{-1}=I +\sum_{d=1}^{+\infty} (-1)^d P^d$$
converges absolutely and we can bound the entries of $Q=(I+P)^{-1}$, $q=\left(q_{il}\right)$, 
by 
$$\split &\left| q_{il} \right| \ \leq \ {1 \over \epsilon^2}  {1 \over \epsilon^{3/2} \sqrt{\delta mn}}
= {1 \over \epsilon^{7/2} \sqrt{\delta mn}}
\quad \text{if} \quad i \ne l \\ &\qquad \quad \text{and} \\
&\left|q_{ii} -1 \right| \ \leq \   {1 \over \epsilon^{7/2} \sqrt{\delta mn}}. \endsplit$$
On the other hand, $Q$ is the matrix of covariances of functions $s_1, \ldots, s_m; t_1, \ldots, t_n$, 
so we have 
$$\aligned &\left|\EE \left(s_j^2; \ \tilde{\psi} \right)-1\right|, \ \left|\EE\left(t_k^2;\ \tilde{\psi} \right)-1 \right| \ \leq \ 
 {1 \over \epsilon^{7/2} \sqrt{\delta mn}} \quad \text{for all} \quad j, k, \\
&\left|\EE \left(s_{j_1} s_{j_2};\ \tilde{\psi} \right)\right| \ \leq \ {1 \over \epsilon^{7/2} \sqrt{\delta mn}}  \quad \text{if} \quad j_1 \ne j_2, \\ &\left|\EE \left(t_{k_1} t_{k_2};\ \tilde{\psi} \right)\right|
 \ \leq \  {1 \over \epsilon^{7/2} \sqrt{\delta mn}} \quad 
\text{if} \quad   k_1 \ne k_2 \\ &\qquad \qquad  \text{and} \\
&\left|\EE\left(s_j t_k; \ \tilde{\psi} \right) \right| \ \leq \  {1 \over \epsilon^{7/2} \sqrt{\delta mn}}
\quad \text{for all} \quad j, k. \endaligned
\tag3.3.8$$

Now we go from $\tilde{\psi}$ back to $\psi$. Since $v$ and $w$ are eigenvectors of 
$\tilde{\psi}$ and since $L_0 =(v, w)^{\bot}$,  for any linear functions 
$\ell_1, \ell_2: {\Bbb R}^{m+n} \longrightarrow {\Bbb R}$ we have 
$$\EE \left(\ell_1 \ell_2; \ \tilde{\psi}|L_0 \right) = \EE\left(\ell_1 \ell_2; \ \tilde{\psi}\right) -
\EE\left(\ell_1 \ell_2;\ \tilde{\psi}|\spa(w) \right)-
\EE \left(\ell_1 \ell_2;\ \tilde{\psi}| \spa(v) \right).$$
On the other hand, since $\psi$ and $\tilde{\psi}$ coincide on $L_0$, we have
$$\EE\left(\ell_1 \ell_2;  \ {\psi}|L_0 \right)=\EE\left(\ell_1 \ell_2;  \ \tilde{\psi}|L_0 \right).$$
Finally, since $v$ is an eigenvector of $\psi$ and $L_0$ is the orthogonal complement to $v$ in 
$L$, we have 
$$\EE \left(\ell_1 \ell_2;\ \psi|L \right) =\EE\left(\ell_1 \ell_2;\ \psi|L_0 \right) + \EE\left(\ell_1 \ell_2; \ \psi|\spa(v) \right).$$
Therefore,
$$\aligned \EE \left(\ell_1 \ell_2; \ \psi|L \right) =&\EE\left(\ell_1 \ell_2;\  \tilde{\psi} \right)-
\EE\left(\ell_1\ell_2;\ \tilde{\psi}|\spa(w) \right) \\ &-
\EE \left(\ell_1 \ell_2;\ \tilde{\psi}| \spa(v) \right)+ \EE\left(\ell_1 \ell_2;\ \psi|\spa(v) \right). \endaligned \tag3.3.9$$
We note that the gradient of function $s_j$ restricted onto $\spa(w)$ is $\sqrt{a_j/2c}$. Since 
$w$ is an eigenvector of $\tilde{\psi}$ with eigenvalue $\epsilon^2/2$, we have 
$$\left(\EE s_{j_1} s_{j_2}; \ \tilde{\psi}|\spa(w) \right)={\sqrt{a_{j_1} a_{j_2}} \over 2 \epsilon^2 c} \ \leq \ 
{1 \over 2 \epsilon^3 \sqrt{\delta mn}}
\quad \text{for all} 
\quad j_1, j_2.$$
Similarly,
$$\left(\EE t_{k_1} t_{k_2};\ \tilde{\psi}| \spa(w)\right) ={\sqrt{b_{k_1} b_{k_2}} \over 2 \epsilon^2 c}
\ \leq  \ {1 \over 2 \epsilon^3 \sqrt{\delta mn}} \quad \text{for all} \quad k_1, k_2
$$
and 
$$\left(\EE s_j t_k; \ \tilde{\psi}| \spa(w) \right) =-{\sqrt{a_j b_k} \over 2 \epsilon^2 c}  \ \geq \ 
-{1 \over 2\epsilon^3 \sqrt{ mn}} \quad \text{for all} \quad j, k.$$
Since $v$ is an eigenvector of $\tilde{\psi}$ with eigenvalue $1-\epsilon^2/2 \geq 1/2$, we obtain 
$$\left(\EE s_{j_1} s_{j_2};\ \tilde{\psi}|\spa(v) \right)=
{\sqrt{a_{j_1} a_{j_2}} \over 4 \left(1-\epsilon^2/2\right) c} \ \leq \ 
{1 \over 2 \epsilon \sqrt{\delta mn}}
\quad \text{for all} 
\quad j_1, j_2.$$
Similarly,
$$\left(\EE t_{k_1} t_{k_2}; \ \tilde{\psi}| \spa(v)\right) ={\sqrt{b_{k_1} b_{k_2}} \over 4 
\left(1-\epsilon^2/2 \right) c}
\ \leq  \ {1 \over 2\epsilon \sqrt{\delta mn}} \quad \text{for all} \quad k_1, k_2
$$
and 
$$\left(\EE s_j t_k; \ \tilde{\psi}| \spa(v) \right) ={\sqrt{a_j b_k} \over 4 \left(1- \epsilon^2/2 \right) c}  \ \leq \ 
{1 \over 2\epsilon \sqrt{mn}} \quad \text{for all} \quad j, k.$$
Since $v$ is an eigenvector of $\psi$ with eigenvalue 1, we get 
$$\left(\EE s_{j_1} s_{j_2}; \ \psi|\spa(v) \right)=
{\sqrt{a_{j_1} a_{j_2}} \over 4  c} \ \leq \ 
{1 \over 4 \epsilon \sqrt{\delta mn}}
\quad \text{for all} 
\quad j_1, j_2.$$
Similarly,
$$\left(\EE t_{k_1} t_{k_2}; \ \psi| \spa(v)\right) ={\sqrt{b_{k_1} b_{k_2}} \over 4 c}
\ \leq  \ {1 \over 4 \epsilon \sqrt{\delta mn}} \quad \text{for all} \quad k_1, k_2
$$
and 
$$\left(\EE s_j t_k; \ \psi| \spa(v) \right) ={\sqrt{a_j b_k} \over 4 c}  \ \leq \ 
{1 \over 4 \epsilon \sqrt{mn}} \quad \text{for all} \quad j, k.$$
Combining (3.3.8) and (3.3.9), we complete the proof.
{\hfill \hfill \hfill} \qed
\enddemo

Now we are ready to prove Theorem 3.2.

\subhead (3.4) Proof of Theorem 3.2 \endsubhead
Let us define
$$\xi_{jk} =\zeta_{jk}^2 + \zeta_{jk} \quad \text{for all} \quad j, k.$$
Hence we have 
$$\split &\alpha  \ \leq \ \xi_{jk} \ \leq \ \beta \quad \text{for all} \quad j, k, 
\quad \text{where} \\ &\alpha =\tau \delta + \tau^2 \delta^2  \quad \text{and} \quad 
\beta = \tau + \tau^2. \endsplit$$
We have 
$${\beta \over \alpha} ={\tau + \tau^2 \over \tau \delta + \tau^2 \delta^2} ={1+\tau \over \delta + \tau \delta^2} \ \leq \ {1 \over \delta^2}. \tag3.4.1$$
Let
$$a_j =\sum_{k=1}^n \xi_{jk} \quad \text{and} \quad b_k=\sum_{j=1}^m \xi_{jk}.$$
In particular, we have 
$$\aligned &a_j \ \leq \ \left(\tau + \tau^2 \right) n \quad \text{for} \quad j=1, \ldots, m \quad \text{and} \\
 &b_k \ \leq \ \left(\tau +\tau^2 \right) m \quad \text{for} \quad k=1, \ldots, n. \endaligned \tag3.4.2$$
We apply Proposition 3.3 to the quadratic form 
$$\psi={1 \over 2} \sum \Sb 1 \leq j \leq m \\ 1 \leq k \leq n \endSb \xi_{jk}
 \left({s_j \over \sqrt{a_j}} +{t_k \over \sqrt{b_k} }\right)^2$$ and the hyperplane 
$L_1 \subset {\Bbb R}^{m+n}$ defined by the equation
$$\sum_{j=1}^m s_j \sqrt{a_j}\ =\ \sum_{j=1}^n t_k \sqrt{b_k}.$$
Let us consider a linear transformation 
$$\left(s_1, \ldots, s_m; \ t_1, \ldots, t_n \right) \longmapsto \left(s_1 \sqrt{a_1}, \ldots ,
s_m \sqrt{a_m}; \ t_1 \sqrt{b_1}, \ldots, t_n \sqrt{b_n} \right)$$
and the hyperplane $L_2 \subset {\Bbb R}^{m+n}$ defined by the equation
$$\sum_{j=1}^m a_j s_j =\sum_{k=1}^n b_k t_k.$$
Then $L_2$ is mapped onto  $L_1$ and the push-forward of the Gaussian probability measure on $L_2$ 
with the density proportional to $e^{-q}$ is the Gaussian probability measure on $L_1$ with the density 
proportional to $e^{-\psi}$. 

We have 
$$\aligned &\EE\left(s_{j_1} s_{j_2}; \ q| L_2 \right)= {1 \over \sqrt{a_{j_1} a_{j_2}}}
 \EE\left(s_{j_1} s_{j_2}; \  
\psi|L_1 \right) \quad \text{for all} \quad j_1,  j_2, \\
&\EE\left(t_{k_1} t_{k_2}; \ q| L_2 \right)= {1 \over \sqrt{b_{k_1} b_{k_2}}} \EE\left(t_{k_1} t_{k_2};\ 
\psi|L_1 \right) \quad \text{for all} \quad k_1,  k_2, \quad \text{and} \\
&\EE\left(s_j t_k;\ q| L_2 \right)= {1 \over \sqrt{a_j b_k}} \EE\left(s_j  t_k;\
\psi|L_1 \right) \quad \text{for all} \quad j, k. \endaligned \tag3.4.3$$
By (3.4.1), we have $\beta/\alpha \leq \delta^{-2}$.
Since by Lemma 3.1, for any hyperplane $L \subset {\Bbb R}^{m+n}$ not containing $u$
we have 
$$\EE \Bigl(\left(s_{j_1} + t_{k_1}\right)\left(s_{j_2}+t_{k_2}\right); \ q|L \Bigr)=
\EE \Bigl(\left(s_{j_1} + t_{k_1}\right)\left(s_{j_2}+t_{k_2}\right);\  q|L_2 \Bigr), $$ the proof follows by 
Proposition 3.3 applied to $L_1$ and $\psi$ and (3.4.1)--(3.4.3).
{\hfill \hfill \hfill} \qed

We will need the following result.
\proclaim{(3.5) Lemma} Let $V$ be Euclidean space and let $q: V \longrightarrow {\Bbb R}$ be a quadratic form such that $\rk q=\dim V-1$. 
Let $v \in V$ be the unit eigenvector of $q$ with the eigenvalue $0$ and let $H=v^{\bot}$ be the 
orthogonal complement of $v$. Then for a unit vector $u \in V$ we have
$$\det q|u^{\bot}=\langle u, v \rangle^2 \det q|H.$$
\endproclaim
\demo{Proof}
This is Lemma 2.3 of \cite{B97b}.
{\hfill \hfill \hfill} \qed
\enddemo

We apply Lemma 3.5 in the following situation. Let $V={\Bbb R}^{m+n}$ and let $q$ be defined by
(1.2.1). Let $L$ be a coordinate hyperplane defined by one of the equations $s_j=0$ or $t_k=0$.
Then 
$$\det q|L = {1 \over m+n} \det q|H.$$
In particular, the value of $\det q|L$ does not depend on the choice of the coordinate hyperplane.

Finally, we need the following result.
\proclaim{(3.6) Lemma} Let 
$q_0: {\Bbb R}^{m+n} \longrightarrow {\Bbb R}$ be the quadratic form defined by the formula
$$q_0(s, t)={1 \over 2} \sum \Sb 1 \leq j \leq m \\ 1 \leq k \leq n \endSb \left(s_j +t_k \right)^2.$$
Then the eigenspaces of $q_0$ are as follows:

The 1-dimensional eigenspace $E_1$ with the eigenvalue 0 spanned by vector 
$$u=\left( \underbrace{1, \ldots, 1}_{\text{$m$ times}}; \ \underbrace{-1, \ldots, -1}_{\text{$n$ times}}
\right);$$
The $(n-1)$-dimensional eigenspace $E_2$ with the eigenvalue $m/2$ consisting of the vectors 
such that 
$$\sum_{k=1}^n t_k =0 \quad \text{and} \quad s_1 =\ldots =s_m =0;$$
The $(m-1)$-dimensional eigenspace $E_3$ with the eigenvalue $n/2$ consisting of the vectors 
such that 
$$\sum_{j=1}^m s_j =0 \quad \text{and} \quad t_1 =\ldots = t_n =0$$
and

\noindent The 1-dimensional eigenspace $E_4$ with the eigenvalue $(m+n)/2$ spanned by vector 
$$v=\left( \underbrace{n, \ldots, n}_{\text{$m$ times}}; \ \underbrace{m, \ldots, m}_{\text{$n$ times}}
\right).$$
\endproclaim
\demo{Proof} Clearly, $E_1$ is the eigenspace with the eigenvalue 0. It is then straightforward to check 
that the gradient of $q_0$ at $x=(s, t)$ equals $mx$ for $x \in E_2$, equals $nx$ for $x \in E_3$ 
and equals $(m+n)x$ for $x\in E_4$.
{\hfill \hfill \hfill} \qed
\enddemo

\head 4. The third degree term \endhead

In this section we prove the following main result.
\proclaim{(4.1) Theorem} Let $u_{jk}$, $j=1, \ldots, m$ and $k=1, \ldots, n$ be Gaussian random 
variables such that 
$$\EE u_{jk}=0 \quad \text{for all} \quad j, k.$$
Suppose further that for some $\theta>0$ 
$$\split &\EE u_{jk}^2 \ \leq \ {\theta \over m+n} \quad \text{for all} \quad j, k \quad \text{and} \\
&\left| \EE u_{j_1 k_1} u_{j_2 k_2} \right| \ \leq \ {\theta \over  (m+n)^2} \quad \text{provided} \quad 
j_1 \ne j_2 \quad \text{and} \quad k_1 \ne k_2. \endsplit$$
Let 
$$U =\sum \Sb 1 \leq j \leq m \\ 1 \leq k \leq n \endSb u_{jk}^3.$$
Then for some constant $\gamma(\theta) >0$ and any $0 < \epsilon \leq 1/2$ we have
$$\left| \EE \exp\{i U \}  - \exp\left\{ - {1 \over 2} \EE U^2 \right\} \right| \ \leq \  \epsilon$$
provided 
$$m + n \ \geq \ \left({1 \over \epsilon}\right)^{\gamma(\theta)}.$$
Besides,
$$\EE U^2 \ \leq \ \gamma(\theta)$$
for some constant $\gamma(\theta) >0$.
Here $i=\sqrt{-1}$.
\endproclaim
We will apply Theorem 4.1 in the following situation. Let $q: {\Bbb R}^{m+n} \longrightarrow {\Bbb R}$ 
be the quadratic form defined by (1.2.1). Let $L \subset {\Bbb R}^{m+n}$ be a hyperplane 
not containing the null-space of $q$. Let us fix the Gaussian probability measure on $L$ with the 
density proportional to $e^{-q}$. We define random variables 
$$u_{jk} = \root 3 \of{{\zeta_{jk} \left(\zeta_{jk}+1 \right) \left(2 \zeta_{jk} +1 \right) \over 6}}
\left(s_j + t_k\right),$$
where $s_1, \ldots, s_m; t_1, \ldots, t_n$ are the coordinates of a point in $L$.
Then we have 
$$U=\sum \Sb 1 \leq j  \leq m \\ 1 \leq k \leq n  \endSb u_{jk}^3 =f(s,t)$$
for $f$ defined by (1.2.2).

\subhead (4.2) The expectation of a product of Gaussian random variables \endsubhead
We will use the famous Wick's formula, see, for example, \cite{Zv97}.
Let $w_1, \ldots, w_l$ be Gaussian random variables such that 
$$\EE w_1 = \ldots = \EE w_l=0.$$ 
Then 
$$\split &\EE w_1 \cdots w_l =0 \quad \text{if} \quad l=2r+1 \quad \text{is odd} \qquad \text{and} \\
&\EE w_1 \cdots w_l = \sum \left(\EE w_{i_1} w_{i_2} \right) \cdots \left(\EE w_{i_{2r-1}} w_{i_{2r}} \right)
\quad \text{if} \quad l=2r \quad \text{is even} \endsplit$$
and the sum is taken over all $(2r)!/r! 2^r$ unordered partitions of the set of
 indices $\{1, \ldots, 2r\}$ into 
$r$ pairwise disjoint unordered pairs $\left\{i_1, i_2 \right\}, \ldots, \left\{i_{2r-1}, i_{2r} \right\}$.
Such a partition is called a {\it matching} of the random variables $w_1, \ldots, w_l$.
We say that $w_i$ and $w_j$ are {\it matched} if they form a pair in the matching.

In particular,
$$\EE w^{2r}={(2r)! \over r! 2^r} \left( \EE w^2 \right)^r. \tag4.2.1$$
We will also use that 
$$\EE \left(w_1^3 w_2^3 \right) = 9 \left(\EE w_1^2 \right) 
\left(\EE w_2^2 \right) \left( \EE w_1 w_2 \right) +
6 \left( \EE w_1 w_2\right)^3 \tag4.2.2$$
and later in Section 5  that
$$\aligned\cov \left(w_1^4, w_2^4 \right) =&\EE \left(w_1^4 w_2^4 \right) -\left(\EE w_1^4 \right)
\left(\EE w_2^4 \right)\\=&9 \left(\EE w_1^2\right)^2 \left(\EE w_2^2 \right)^2+
72 \left(\EE w_1 w_2 \right)^2 \left(\EE w_1^2 \right) \left(\EE w_2^2 \right)\\
&\qquad +24 \left(\EE w_1 w_2\right)^4 -9\left(\EE w_1^2 \right)^2 \left(\EE w_2^2 \right)^2\\ =
&72 \left(\EE w_1 w_2 \right)^2 \left(\EE w_1^2 \right) \left(\EE w_2^2 \right)+24 \left(\EE w_1 w_2\right)^4
\endaligned \tag4.2.3$$

\bigskip
All implied constants in the ``$O$'' notation in this section are absolute.
\bigskip

\subhead (4.3) Auxiliary random variables $v_{jk}$ \endsubhead 
For the Gaussian random variables $\left\{u_{jk}\right\}$ of Theorem 4.1, 
let us define Gaussian random variables $\left\{ v_{jk} \right\}$, where $j=1, \ldots, m$ and 
$k=1, \ldots, n$ such that 
$$\split & \EE v_{jk} =0 \quad \text{for all} \quad j, k \quad \text{and} \\
& \EE \left(v_{j_1 k_1} v_{j_2 k_2}\right) =\EE \left(u_{j_1 k_1}^3 u_{j_2 k_2}^3\right) \quad \text{for all} \quad 
j_1, j_2 \quad \text{and} \quad k_1, k_2. \endsplit$$
We say that the random variables $u_{j_1k_1}$ and $u_{j_2k_2}$ in Theorem 4.1 are 
{\it weakly correlated} if $j_1 \ne j_2$ and $k_1\ne k_2$ and that they are {\it strongly correlated}
otherwise.  Similarly, we say that $v_{j_1 k_1}$ and $v_{j_2 k_2}$ are {\it weakly correlated}
if $j_1 \ne j_2$ and $k_1 \ne k_2$ and are {\it strongly correlated} otherwise.

By (4.2.2), 
$$\aligned & \EE u_{j_1 k_1}^3 u_{j_2 k_2}^3 = \EE v_{j_1k_1} v_{j_2 k_2} \\ &\qquad =\cases O\left(\dsize{\theta^3 \over (m+n)^4} \right) &\text{if\ } u_{j_1 k_1}, u_{j_2 k_2}
\quad \text{are weakly correlated} \\ O\left(\dsize{ \theta^3 \over (m+n)^3}\right) &\text{if\ } u_{j_1 k_1}, u_{j_2 k_2} 
\quad \text{are strongly correlated.} \endcases \endaligned \tag4.3.1$$
Since the number of weakly correlated pairs is $O\left(m^2 n^2 \right)$ while the 
number of strongly correlated pairs is $O\left(m^2 n + n^2 m \right)$, we obtain that 
$$\EE V^2 =\EE U^2 =O\left(\theta^3\right). \tag4.3.2$$

\subhead (4.4) Representation of monomials by graphs \endsubhead
Let $x_{jk}$, $j=1, \ldots, m$, $k=1, \ldots, n$, be formal commuting variables. We interpret a monomial
in $x_{jk}$ combinatorially, as a weighted graph. Let $K_{m,n}$ be the complete bipartite graph 
with $m+n$ vertices and $mn$ edges $(j, k)$ for $j=1, \ldots, m$ and $k=1, \ldots, n$. 
A {\it weighted graph} $G$ is a set of edges $(j, k)$ of $K_{m, n}$ with positive integer weights
$\alpha_{jk}$ on it. With $G$, we associate a monomial
$$t_G(x)=\prod_{(j, k) \in G} x_{jk}^{\alpha_{jk}}.$$
The {\it weight} $\sum_{e \in G} \alpha_e$ of $G$ is the degree of the monomial. We observe that 
for any $p$ there are not more than $r^{O(r)}(m+n)^p$ distinct weighted graphs of weight $2r$ and 
$p$ vertices. We note that pairs of variables $u_{j_1k_1}$ and $u_{j_2k_2}$ and pairs 
of variables $v_{j_1k_1}$ and $v_{j_2 k_2}$ corresponding to the edges $(j_1, k_1)$ and 
$(j_2, k_2)$ in different connected components are weakly correlated.

\proclaim{(4.5) Lemma} 
Given a graph $G$ of weight $2r$, $r>1$, let us represent it as a vertex-disjoint union 
$$G=G_0 \cup G_1,$$
where $G_0$ consists of $s$ isolated edges of weight 1 each and $G_1$ is a graph with no isolated edges of weight 1 (we may have 
$s=0$ and no $G_0$). 

Then
\roster
\item We have
$$\split \left| \EE t_G\left(u_{jk}^3\right)\right| \ \leq \  &{r^{O(r)} \theta^{3r} \over (m+n)^{3r+s/2}} \qquad  \text{and} \\
\left| \EE t_G\left(v_{jk}\right)\right| \ \leq \ &{r^{O(r)}\theta^{3r} \over (m+n)^{3r +s/2}}. \endsplit$$
Additionally, if  $s$ is odd, then 
 $$\left| \EE t_G\left(u_{jk}^3\right)\right|, \  \left| \EE t_G\left(v_{jk}\right)\right| \ \leq \ 
 {r^{O(r)} \theta^{3r} \over (m+n)^{3r+(s+1)/2}}.$$
\item If $s$ is even and $G_1$ has $r-s/2$ connected components, each consisting of a pair of 
edges of weight 1 each sharing one common vertex, then the number of vertices of $G$ is
$$3r + {s \over 2}.$$
In all other cases,  $G$ has strictly fewer than $3r+s/2$ vertices.
\item Suppose that $s$ is even and that $G_1$ has $r-s/2$ connected components, each consisting of 
a pair of 
edges of weight 1 each sharing one common vertex. Then 
$$\left| \EE t_G\left(u_{jk}^3\right)-\EE t_G\left(v_{jk}\right)\right| \ \leq \ {r^{O(r)} \theta^{3r}
\over (m+n)^{3r+s/2+1}}.$$
\endroster
\endproclaim
\demo{Proof} To prove Part (1), we use Wick's formula of Section 4.2. Since for each isolated 
edge $(j_1, k_1) \in G_0$, at least one of the three copies of the random variable $u_{j_1k_1}$ has to be matched with a copy of the variable $u_{j_2 k_2}$ indexed by an edge $(j_2, k_2)$ in a different connected component, we 
conclude that each matching of the multiset 
$$\Bigl\{ u_{jk}, u_{jk}, u_{jk}: \ (j, k) \in G \Bigr\} \tag4.5.1$$
contains at least $s/2$ weakly correlated pairs of variables and hence 
$$\left|\EE t_G\left(u_{jk}^3 \right)\right| \ \leq \ r^{O(r)} \left({\theta \over (m+n)^2} \right)^{s/2} 
\left({\theta \over m+n}\right)^{3r-s/2}.$$
Moreover, if $s$ is odd, the number of weakly correlated pairs in every matching is at least 
$(s+1)/2$ and hence 
$$\left|\EE t_G\left(u_{jk}^3 \right)\right| \ \leq \ r^{O(r)} \left({\theta \over (m+n)^2} \right)^{(s+1)/2} 
\left({\theta \over m+n}\right)^{3r-(s+1)/2}.$$
Similarly, since for each isolated edge $(j_1, k_1) \in G_0$, variable $v_{j_1k_1}$ has to be matched 
with a variable $v_{j_2 k_2}$ indexed by an edge $(j_2, k_2)$ in a different connected component, 
we conclude that each matching of the set 
$$\Bigl\{v_{jk}: \ (j, k) \in G\Bigr\} \tag4.5.2$$
contains at least $s/2$ weakly correlated pairs of variables and hence by (4.3.1) 
$$\left|\EE t_G\left(v_{jk} \right)\right| \ \leq \ r^{O(r)} \left({\theta^3 \over (m+n)^4} \right)^{s/2} 
\left({\theta^3 \over (m+n)^3}\right)^{r-s/2}.$$
Moreover, if $s$ is odd, the number of weakly correlated pairs in every matching is at least $(s+1)/2$ 
and hence
$$\left|\EE t_G\left(v_{jk} \right)\right| \ \leq \ r^{O(r)} \left({\theta^3 \over (m+n)^4} \right)^{(s+1)/2} 
\left({\theta^3 \over (m+n)^3}\right)^{r-(s+1)/2},$$
which concludes the proof of Part (1).

To prove Part (2), we note that a connected weighted graph $G$ of weight $e$ contains a spanning tree with at most $e$
edges and hence has at most $e+1$ vertices. In particular, a connected graph $G$ of weight $e$ 
contains fewer than $3e/2$ vertices unless $G$ is an isolated edge of weight 1 or a pair of edges 
of weight 1 each, sharing one common vertex. Therefore, $G$ has at most 
$$ 2s + {3 \over 2} (2r-s)=3r+{s \over 2}$$
vertices and strictly fewer vertices,
unless $s$ is even and the connected components of $G_1$ are pairs of edges of weight 
1 each sharing one common vertex.

To prove Part (3), let us define $\Sigma_u(G)$ as the sum in the Wick's formula over all matchings 
of the multiset (4.5.1) of the following structure: we split the edges of $G$ into $r$ pairs, pairing each 
isolated edge with another isolated edge and pairing each edge in a 2-edge connected 
component of $G$ with the remaining edge in the same connected component. We then match
every variable of the multiset (4.5.1) with a variable indexed by an edge in the same pair.
Reasoning as in Part (1), we conclude that 
$$\left| \EE t_G\left(u_{jk}^3\right)- \Sigma_u(G)\right| \ \leq \ {r^{O(r)} \theta^{3r} \over (m+n)^{3r+s/2+1}}
.$$
Similarly, let us define $\Sigma_v(G)$ as the sum in the Wick's formula over all matchings 
of the set (4.5.2) of the following structure: we split the edges of $G$ into $r$ pairs as above
and match every variable in the set (4.5.2) with the variable indexed by the remaining edge of 
the pair. Then 
$$\left| \EE t_G\left(v_{jk}\right)- \Sigma_v(G)\right| \ \leq \ {r^{O(r)} \theta^{3r} \over (m+n)^{3r+s/2+1}}
.$$
Since 
$$\Sigma_u(G)=\Sigma_v(G),$$
the proof of Part (3) follows.
{\hfill \hfill \hfill} \qed
\enddemo

\proclaim{(4.6) Lemma} Let $u_{jk}$ be random variables as in Theorem 4.1 and let 
$v_{jk}$ be the auxiliary Gaussian random variables as in Section 4.3.
Let
$$U=\sum \Sb 1 \leq j \leq m \\ 1 \leq k \leq n \endSb u_{jk}^3 \quad \text{and} \quad 
V=\sum \Sb 1 \leq j \leq m \\ 1 \leq k \leq n \endSb v_{jk}.$$
Then 
 for every integer $r>1$ we have 
$$\left| \EE U^{2r} - \EE V^{2r} \right| \ \leq \ {r^{O(r)} \theta^{3r} \over m+n}.$$
\endproclaim
\demo{Proof}
We can write 
$$\split \EE U^{2r} =&\sum_G a_G \EE t_G\left(u_{jk}^3 \right) \quad \text{and} \\
\EE V^{2r}=&\sum_G a_G \EE t_G\left(v_{jk} \right), \endsplit $$
where the sum is taken over all weighted graphs $G$ of the total weight $2r$ and 
$$1 \ \leq \ a_G \ \leq \ (2r)!.$$
Let $\GG_{2r}$ be the set of weighted graphs $G$ whose connected components 
consist of an even 
number $s$ of isolated edges and $r-s/2$ pairs of edges of weight 1 sharing one common vertex.
Since there are not more than $r^{O(r)} (m+n)^p$ distinct weighted graphs with $p$ vertices, 
by Parts (1) and (2) of Lemma 4.5, we conclude that 
$$\split \left|  \EE U^{2r} -\sum_{G \in \GG_{2r}} a_G \EE t_G\left(u_{jk}^3 \right) \right| \ \leq \
&{r^{O(r)} \theta^{3r} \over m+n}  \quad 
\text{and} \\
\left|  \EE V^{2r} -\sum_{G \in \GG_{2r}} a_G \EE t_G\left(v_{jk} \right) \right| \ \leq \
&{r^{O(r)} \theta^{3r} \over m+n}. \endsplit$$
The proof now follows by Parts (1) and (3) of Lemma 4.5.
{\hfill \hfill \hfill} \qed
\enddemo

\subhead (4.7) Proof of Theorem 4.1 \endsubhead 
Let $U$ and $V$ be the random variables as in Lemma 4.6. 
We use the standard estimate 
$$\left| e^{ix} - \sum_{s=0}^{2r-1} {(ix)^s \over s!} \right| \ \leq \ {x^{2r} \over (2r)!} \quad \text{for} \quad 
x \in {\Bbb R},$$
from  which it follows that 
$$\left| \EE e^{iU} -\EE e^{iV} \right| \ \leq \ {\EE U^{2r} \over (2r)!} + {\EE V^{2r} \over (2r)!} 
+\sum_{s=0}^{2r-1} {\left| \EE U^s - \EE V^s \right| \over s!}.\tag4.7.1$$
By (4.3.2), we have 
$$\EE U^2=\EE V^2 =O\left(\theta^3\right)$$ and hence 
$$\EE U^{2r}, \ \EE V^{2r} \ \leq \ {(2r)! \over 2^r r!} 2^{O(r)} \theta^{3r}.$$
Therefore, one can choose an integer $r$ such that 
$$r \ln r \ \leq \ \gamma_1(\theta) \ln {1 \over \epsilon} \quad \text{for some constant} \quad 
\gamma_1(\theta)>0$$
so that 
$${\EE U^{2r} +\EE V^{2r} \over (2r)!} \ \leq \ {\epsilon \over 2}.$$
By Lemma 4.6, as long as 
$$m + n \ \geq \ \left({1 \over \epsilon}\right)^{\gamma(\theta)} \quad \text{for some constant} \quad \gamma(\theta)>0$$
we have 
$$\left| \EE U^{2k} - \EE V^{2k} \right| \leq {\epsilon \over 6} \quad \text{for} \quad k=0, 1, \ldots, r.$$ 
We note that by symmetry 
$$\EE V^s =\EE U^s =0 \quad \text{if} \quad s \quad \text{is odd}.$$
Since $V$ is Gaussian, we have 
$$\EE e^{iV}=\exp\left\{-{1 \over 2} \EE V^2 \right\},$$
and the proof follows by (4.7.1).
{\hfill \hfill \hfill} \qed

\head 5. The fourth degree term \endhead 

In this section we prove the following main result.
\proclaim{(5.1) Theorem} Let $w_{jk}$, $j=1, \ldots, m$ and $k=1, \ldots, n$ be Gaussian random 
variables such that 
$$\EE w_{jk}=0 \quad \text{for all} \quad j, k.$$
Suppose further that for some $\theta >0$ we have 
$$\split &\EE w_{jk}^2 \ \leq \ {\theta \over m+n} \quad \text{for all} \quad j, k \quad \text{and} \\
&\left| \EE w_{j_1 k_1} w_{j_2 k_2} \right| \ \leq \ {\theta \over (m+n)^2} \quad 
\text{provided} \quad j_1 \ne j_1 \quad \text{and} \quad k_1 \ne k_2. \endsplit $$
Let
$$W=\sum \Sb 1 \leq j \leq m \\ 1 \leq k \leq n \endSb w_{jk}^4.$$
Then for some absolute constant $\gamma>0$ we have
\roster
\item $$\EE W \ \leq \ \gamma \theta^2; $$
\item $$\var W  \ \leq \ {\gamma \theta^4 \over m+n};$$
\item $$\PP\Bigl\{W > \gamma \theta^2 +1 \Bigr\} \ \leq \ \exp\left\{ -(m+n)^{1/5} \right\}$$
provided $m + n \geq \gamma_1(\theta)$ for some $\gamma_1(\theta)>0$.
\endroster 
\endproclaim

We will apply Theorem 5.1 in the following situation. Let $q: {\Bbb R}^{m+n} \longrightarrow {\Bbb R}$
be the quadratic form defined by (1.2.1) and let $L \subset {\Bbb R}^{m+n}$ be a hyperplane 
not containing the null-space of $q$. Let us fix the Gaussian probability measure in $L$ with the 
density proportional to $e^{-q}$. We define random variables 
$$w_{jk}=\root 4\of{{ \zeta_{jk} \left(\zeta_{jk} +1 \right) \left(6 \zeta_{jk}^2 +6\zeta_{jk}+1 \right)
\over 24}}\left(s_j +t_k\right),$$
where $s_1, \ldots, s_m; t_1, \ldots, t_n$ are the coordinates of a point in $L$.
Then we have 
$$W=\sum \Sb 1 \leq k \leq m \\ 1 \leq j \leq n \endSb w_{jk}^4 =h(s,t)$$
for $h$ defined by (1.2.2).

While the proof of Parts (1) and (2) is a straightforward computation, to prove Part (3) we need 
reverse H\"older inequalities for polynomials with respect to the Gaussian measure.

\proclaim{(5.2) Lemma} Let $p$ be a polynomial of degree $d$ in random Gaussian
variables $w_1, \ldots, w_l$.

Then for $r > 2$ we have 
$$\left(\EE |p|^r\right)^{1/r} \ \leq \ r^{d/2} \left(\EE p^2 \right)^{1/2}.$$
\endproclaim
\demo{Proof} This is Corollary 5 of \cite{Du87}.
\enddemo
{\hfill \hfill \hfill} \qed

\subhead (5.3) Proof of Theorem 5.1 \endsubhead
Using formula (4.2.1), we get
$$\EE w_{jk}^4 =3 \left( \EE w_{jk}^2 \right)^2 \ \leq \ {3 \theta^2 \over (m+n)^2}$$
and hence
$$\EE W =\sum \Sb 1 \leq j \leq m \\ 1 \leq k \leq n \endSb \EE w_{jk}^4 
\ \leq \ (mn){3 \theta^2 \over (m+n)^2} \ \leq \ 3 \theta^2,$$
which proves Part (1).

To prove Part (2), we note that
$$\var W =\sum \Sb 1 \leq j_1, j_2 \leq m \\ 1 \leq  k_1, k_2 \leq n \endSb 
\cov \left( w_{j_1 k_1}^4, w_{j_2 k_2}^4 \right).$$
Using (4.2.3), we get
$$ \cov \left(w_{j_1 k_1}^4, w_{j_2 k_2}^4 \right)=72 \left(\EE w_{j_1 k_1} w_{j_2 k_2} \right)^2
 \left(\EE w_{j_1 k_1}^2 \right) \left(\EE w_{j_2 k_2}^2 \right) +24 \left(\EE w_{j_1 k_1} w_{j_2 k_2} \right)^4. $$
 Hence
 $$\cov \left(w_{j_1 k_1}^4, w_{j_2 k_2}^4 \right) \ \leq \ {96\theta^4 \over (m+n)^4}
 \quad \text{for all} \quad j_1, j_2 \quad \text{and} \quad k_1, k_2.$$
 Additionally,
 $$\split \cov \left(w_{j_1 k_1}^4, w_{j_2 k_2}^4 \right) \ \leq \ &{72 \theta^4 \over (m+n)^6}+
 {24 \theta^4 \over (m+n)^8} \\ \leq \ &{96 \theta^4 \over (m+n)^6}
 \quad \text{provided} \quad j_1 \ne j_2 \quad \text{and} \quad k_1 \ne k_2. \endsplit$$
 Summarizing,
 $$\var W \ \leq \ m^2 n^2 {96 \theta^4 \over (m+n)^6} + \left(mn^2 +nm^2\right)
 {96 \theta^4 \over (m+n)^4} \ \leq \ {\gamma \theta^4 \over m+n},$$
 which proves Part (2).
 
 To prove Part (3), we apply Lemma 5.2 with $d=4$ to $p(W)=W -\EE W$. From Part (2), we get 
$$\EE \left| W-\EE W \right|^r \ \leq \ r^{2r} \left(\var W\right)^{r/2} \ \leq \ 
r^{2r} \left({\gamma \theta^4 \over m+n} \right)^{r/2}. $$
Let us choose
$$r=(m+n)^{1/5}$$
for sufficiently large $m+n$.
Then
$$\split r^{2r} \left({\gamma \theta^4 \over m+n}\right)^{r/2}=
&\exp\left\{ 2r \ln r + {r \over 2} \ln \left(\gamma \theta^4\right) - {r \over 2}\ln (m+n) \right\} \\
=&\exp\left\{-{1 \over 10} (m+n)^{1/5} \ln (m+n) +{\ln \left(\gamma \theta^4\right) \over 2} (m+n)^{1/5} 
\right\} \\
\leq & \exp\left\{- {(m+n)^{1/5} \over 10} \right\}\endsplit$$
provided $\ln (m+n) > 5 \ln \left(\gamma \theta^4 \right) +10$

Hence if $m+n$ is sufficiently large, we have 
$$\EE \left| W - \EE W \right|^r \ \leq \ \exp\left\{ -{(m+n)^{1/5} \over 10}\right\}.$$
The proof follows by Part (1) and Markov's inequality.
 {\hfill \hfill \hfill} \qed

\head 6.  Computing the integral over a neighborhood of the origin \endhead

We consider the integral 
$$\int_{\Pi_0} F(s, t) \ ds dt$$  of Corollary 2.2. Recall that $\Pi_0$ is the facet 
of the parallelepiped
$$\Pi=\Bigl\{ \left(s_1, \ldots, s_m; t_1, \ldots, t_n \right): \quad -\pi \leq s_j, t_k \leq \pi \quad
\text{for all} \quad j, k \Bigr\}$$
defined by the equation $t_n=0$ and that 
$$F(s, t)=\exp\left\{-i \sum_{j=1}^m r_j s_j -i \sum_{k=1}^n c_k t_k \right\}
\prod\Sb 1 \leq j \leq m \\ 1 \leq k \leq n \endSb {1 \over 1+\zeta_{jk} -\zeta_{jk} e^{i(s_j +t_k)}}.$$
In this section, we prove the following main result. 

\proclaim{(6.1) Theorem} Let us fix  a number $0 < \delta < 1$.
Suppose that $m \geq \delta n$, $n \geq \delta m$ and that 
$$ \delta \tau \ \leq \ \zeta_{jk} \ \leq \ \tau \quad \text{for all} \quad j, k$$
and some $\tau >\delta$. 

Let $q: {\Bbb R}^{m+n} \longrightarrow {\Bbb R}$ be the quadratic form defined by (1.2.1)
and \break let $f, h: {\Bbb R}^{m+n} \longrightarrow {\Bbb R}$ be the polynomials defined by (1.2.2).
Let us define a neighborhood $\UU \subset \Pi_0$ by
$$\UU=\left\{(s, t) \in \Pi_0: \quad |s_j|, |t_k| \ \leq \ {\ln (m+n) \over \tau \sqrt{m+n}} \quad \text{for all} \quad
j, k \right\}.$$  
Let us identify the hyperplane $t_n=0$ containing $\Pi_0$ with ${\Bbb R}^{m+n-1}$, let 
$$\Xi =\int_{{\Bbb R}^{m+n-1}} e^{-q} \ ds dt$$
and let 
us consider the Gaussian probability measure in ${\Bbb R}^{m+n-1}$ with the density
$\Xi^{-1} e^{-q}$. Let 
$$\mu = \EE f^2 \quad \text{and} \quad  \nu =\EE h.$$
Then 
\roster
\item  $$\Xi \ \geq \  {(2 \pi)^{(m+n-1)/2} \over m^{(n-1)/2} n^{(m-1)/2} (\tau + \tau^2)^{(m+n-1)/2}}.$$
\item $$0 \ \leq \ \mu, \ \nu \ \leq \gamma(\delta)$$ 
for some constant $\gamma(\delta)>0$.
\item For any $0 <\epsilon \leq 1/2$
$$ \left| \int_{\UU} F(s, t) \ ds dt  -\exp\left\{-{\mu \over 2} +\nu \right\} \Xi  \right| \ \leq \ 
\epsilon \thinspace \Xi$$
provided 
$$m+n \ \geq \ \left({1 \over \epsilon}\right)^{\gamma(\delta)}$$
for some constant $\gamma(\delta) >0$.
\item For any $0 \leq \epsilon \leq 1/2$
$$ \left| \int_{\UU} |F(s, t)| \ ds dt  - \exp\{\nu\}  \Xi  \right| \ \leq \ 
\epsilon\thinspace \Xi$$
provided 
$$m+n \ \geq \ \left({1 \over \epsilon}\right)^{\gamma(\delta)}$$
for some constant $\gamma(\delta) >0$.
\endroster
\endproclaim 
\demo{Proof} All implied constants in the ``$O$'' and ``$\Omega$'' notations below may depend only on the parameter $\delta$.

Let 
$$q_0(s, t)={1 \over 2} \sum \Sb 1 \leq j \leq m \\ 1 \leq k \leq n \endSb \left(s_j +t_k \right)^2$$
as in Lemma 3.6.
Then
$$q(s, t) \ \leq \ (\tau+ \tau^2) q_0(s,t)$$
and, therefore,
$$\split \int_{{\Bbb R}^{m+n-1}} \exp\{-q(s,t) \}  \ ds dt \ \geq\ &\int_{{\Bbb R}^{m+n-1}} 
\exp\left\{-(\tau+ \tau^2) q_0(s, t) \right\} \ ds dt \\=
&{\pi^{(m+n-1)/2} \over (\tau^2 +\tau)^{(m+n-1)/2} \sqrt{\det q_0 | {\Bbb R}^{m+n-1}}}, \endsplit$$
where $q_0 | {\Bbb R}^{m+n-1}$ is the restriction of the form $q_0$ onto the coordinate 
hyperplane $t_n=0$ in ${\Bbb R}^{m+n}$.

Let $H \subset {\Bbb R}^{m+n}$ be the orthogonal complement complement in ${\Bbb R}^{m+n}$ 
to the kernel of $q_0$, that is, the hyperplane defined by the equation:
$$s_1 +\ldots s_m =t_1 + \ldots + t_n  \quad \text{for} \quad (s, t)=\left(s_1, \ldots, s_m; t_1, \ldots, t_m\right).$$ Then, from the eigenspace decomposition of Lemma 3.6 it follows that 
$$\det q_0 | H = \left({m \over 2}\right)^{n-1} \left({n \over 2}\right)^{m-1} \left( {m+n \over 2}\right).$$
On the other hand, by Lemma 3.5,
$$\det q_0| {\Bbb R}^{m+n-1} ={1 \over m+n} \det q_0 |H$$
and the proof of Part (1) follows.

Let us consider the coordinate functions $s_j, t_k$ as random variables on the 
space ${\Bbb R}^{m+n -1}$ with the Gaussian probability measure with the density $\Xi^{-1} e^{-q}$. 
From Theorem 3.2, 
$$\aligned & \EE \left(s_j +t_k \right)^2 =O\left( {1 \over \tau^2(m+n)}\right)\quad 
\text{for all} \quad j, k\\ & \qquad \qquad \text{and} \\
& \left| \EE \left(s_{j_1} + t_{k_1} \right) \left(s_{j_2} +t_{k_2} \right) \right| =
O\left({1 \over \tau^2 mn}\right) \\ &\qquad \qquad \text{if} \quad j_1 \ne j_2 \quad 
\text{and} \quad k_1 \ne k_2 \endaligned \tag6.1.1 $$ 
Let
$$\aligned u_{jk} =&\root 3 \of{{\zeta_{jk}\left(\zeta_{jk}+1\right)\left(2\zeta_{jk}+1\right) \over 6}}
\left(s_j +t_k \right) \\ &\qquad \text{and} \\
w_{jk}=&\root 4 \of{{\zeta_{jk}\left(\zeta_{jk}+1\right)\left(6 \zeta_{jk}^2 + 6 \zeta_{jk} +1 \right)
\over 24}}\left(s_j + t_k \right) \quad \text{for all} \quad j, k. \endaligned $$
Then $u_{jk}$ and $w_{jk}$ are Gaussian random variables such that
$$\aligned \EE u_{jk}^2, \ \EE w_{jk}^2 = &O\left({1 \over m+n}\right) \quad \text{for all} \quad j, k \\
& \text{and} \\
\left| \EE u_{j_1 k_1} u_{j_2 k_2} \right|, \ \left|\EE w_{j_1 k_1} w_{j_2 k_2} \right| =
& O\left({1 \over mn}\right) \\ &\text{if} \quad j_1 \ne j_2 \quad \text{and} \quad k_1 \ne k_2.
\endaligned \tag6.1.2$$
We observe that 
$$f=\sum \Sb 1\leq j \leq m \\ 1 \leq k \leq n \endSb u_{jk}^3 \quad \text{and} \quad 
h=\sum \Sb 1 \leq j \leq m \\ 1 \leq k \leq n \endSb w_{jk}^4.$$
Therefore, the bound
$$\mu =\EE f^2 =O(1)$$ follows by Theorem 4.1 while the bound 
$$\nu=\EE h=O(1)$$ follows by 
Part (1) of Theorem 5.1. This concludes the proof of Part (2).

Since $s_j + t_k$ is a Gaussian random variable, by the first inequality of (6.1.1), we conclude that 
$$\PP\left\{ \left| s_j + t_k \right| \ \geq \ { \ln (m+n) \over 4\tau \sqrt{m+n}} \right\} 
\ \leq \ \exp\Bigl\{ -\Omega\left( \ln^2(m+n)\right) \Bigr\}.$$
Note that the inequalities hold for $k=n$ with $t_n=0$ as well, and hence
$$\int_{{\Bbb R}^{m+n-1}\setminus \UU} e^{-q} \ ds dt
 \ \leq \ \exp\Bigl\{-\Omega\left( \ln^2(m+n)\right) \Bigr\} \Xi, \tag6.1.3$$
where $\UU$ is the neighborhood defined in the theorem.
Therefore,
$$\left| \int_{{\Bbb R}^{m+n-1}\setminus \UU} \exp\bigl\{-q -i f \bigr\} \ ds dt \right|
 \ \leq \ \exp\Bigl\{-\Omega\left( \ln^2(m+n)\right) \Bigr\} \Xi.$$
Since $f=\sum_{jk} u_{jk}^3$ and Gaussian random variables $u_{jk}$ satisfy (6.1.2), from Theorem 4.1 we conclude that for any
$0 \leq \epsilon \leq 1/2$ we have
$$\aligned &\left| \int_{{\Bbb R}^{m+n-1}} \exp\{-q-if \} \ ds dt - \exp\left\{-{\mu \over 2}\right\} \right| \ \leq \ \epsilon \thinspace \Xi \\
&\qquad \text{provided} \quad m+n > \left( {1 \over \epsilon} \right)^{O(1)}. \endaligned$$
Therefore, for any $0 \leq \epsilon \leq 1/2$ we have
$$\aligned &\left| \int_{\UU} \exp\{-q-if \} \ ds dt - \exp\left\{-{\mu \over 2}\right\} \right| \ \leq \ \epsilon \thinspace \Xi \\
&\qquad \text{provided} \quad m+n > \left( {1 \over \epsilon} \right)^{O(1)}. \endaligned \tag6.1.4$$

Since $h=\sum_{jk} w_{jk}^4$ and Gaussian random variables $w_{jk}$ satisfy (6.1.2), by 
Part (2) of Theorem 5.1, we have 
$$\var h =\EE (h-\nu)^2 = O\left({1 \over m+n}\right).$$
Applying Chebyshev's inequality, we conclude that for any $\epsilon >0$
$$\int \limits \Sb (s, t) \in \UU:\\ |h(s, t) -\nu| \geq \epsilon \endSb e^{-q}\ ds dt =
O \left({1 \over \epsilon^2 (m+n)}\right)
\Xi. \tag6.1.5$$ 
By Part (3) of Theorem 5.1, for some constant $\gamma(\delta) >0$ we have
$$\int \limits \Sb (s, t) \in \UU:\\ h(s, t) \geq \gamma(\delta) \endSb e^{-q}\ ds dt =
O \left(\exp\bigl\{ -(m+n)^{1/5} \bigr\}\right) \Xi. \tag6.1.6$$
In addition,
$$h = O\left( \ln^4(m+n) \right) \quad \text{in} \quad \UU. \tag6.1.7$$
In view of (6.1.3) and Part (2) of the theorem, (6.1.5)--(6.1.7) imply for any $0 \leq \epsilon \leq 1/2$ we have
$$\aligned &\left| \int_{\UU} \exp\{-q+h\} \ ds dt - \exp\left\{\nu \right\} \right| \ \leq \ \epsilon \thinspace \Xi \\
&\qquad \text{provided} \quad m+n > \left( {1 \over \epsilon} \right)^{O(1)}. \endaligned \tag6.1.8$$
Similarly, from (6.1.4)--(6.1.7) we deduce that 
$$\aligned &\left| \int_{\UU} \exp\{-q-if+h\} \ ds dt - \exp\left\{-{\mu \over 2}+\nu \right\} \right| \ \leq \ \epsilon \thinspace \Xi \\
&\qquad \text{provided} \quad m+n > \left( {1 \over \epsilon} \right)^{O(1)}. \endaligned \tag6.1.9$$
From the Taylor series expansion, cf. (2.2.1), we write 
$$\split F(s, t) =&\exp\bigl\{-q(s,t) -i f(s,t) +h(s,t) + \rho(s,t) \bigr\}  \quad \text{where} \\
&\quad |\rho(s,t)| =O\left( {\ln^5 (m+n) \over \sqrt{m+n}}\right) \quad \text{for} \quad (s, t) \in \UU.
\endsplit$$
Therefore, using (6.1.8) and Part (2) of the theorem we conclude that 
$$\aligned &\left|  \int_{\UU} F \ ds dt  - 
\int_{\UU} \exp\left\{-q-if +h\right\} \ ds dt \right| =O\left({\ln^5 (m+n) \over \sqrt{m+n}}\right) \Xi \\
&\qquad \qquad \qquad  \text{and} \\ 
&\left|  \int_{\UU} |F| \ ds dt  - 
\int_{\UU} \exp\left\{-q+h\right\} \ ds dt \right| = O\left({\ln^5 (m+n) \over \sqrt{m+n}}\right)\Xi. 
\endaligned $$
We complete the proof of Parts (3) and (4) from (6.1.8) and (6.1.9).
{\hfill \hfill \hfill} \qed 
\enddemo

\head 7. Bounding the integral outside of a neighborhood of the origin \endhead 

We consider the integral representation of Corollary 2.2.
In this section we prove that the integral of $F(s,t)$ outside of the neighborhood $\UU$ of the 
origin is asymptotically negligible (note that by Theorem 6.1 the integral
of $F$ and the integral of $|F|$ over $\UU$ have the same order). We prove the following main result. 

\proclaim{(7.1) Theorem} Let us fix a number $0 < \delta < 1$. Suppose that $m \geq \delta n$, 
$n \geq \delta m$ and that 
$$\delta \tau \ \leq \ \zeta_{jk} \ \leq \ \tau \quad \text{for all} \quad j, k$$
and some $\delta \leq \tau \leq  (m+n)^{1/\delta}$.
Let 
$$\UU=\left\{(s, t) \in \Pi_0: \quad |s_j|, |t_k| \ \leq \ {\ln (m+n) \over \tau \sqrt{m+n}} \quad \text{for all} \quad
j, k \right\}.$$  
Then for any $\kappa>0$ 
$$\int_{\Pi_0 \setminus \UU} |F(s, t)| \ ds dt \ \leq \ (m+n)^{-\kappa}  \int_{\Pi_0} |F(s, t)| \ ds dt$$
provided $m+n > \gamma(\delta, \kappa)$ for some $\gamma(\delta, \kappa)>0$.
\endproclaim

We prove Theorem 7.1 it in two steps. First, by a string of combinatorial arguments we show that the integral
$$\int_{\Pi_0 \setminus I} |F(s, t)| \ ds dt$$
is negligible compared to 
$$\int_I |F(s, t)| \ ds dt, \tag7.1.1$$
where $I \subset \Pi_0$ is a larger neighborhood of the origin,
$$I=\Bigl\{(s, t) \in \Pi_0: \quad |s_j|, |t_k| \ \leq \ \epsilon/\tau \quad \text{for all} \quad j, k \Bigr\}$$
and $\epsilon >0$ is any fixed number. This is the only place where we use that $\tau$ is bounded
above by a polynomial in $m+n$.
Then we notice that for a sufficiently small
$\epsilon=\epsilon(\delta)$, the function $|F(s,t)|$ is strictly log-concave on $I$ and we use 
a concentration inequality for strictly log-concave measures to deduce that the integral
$$\int_{I \setminus \UU} |F(s, t)| \ ds dt$$
is negligible compared to (7.1.1).

\subhead (7.2) Metric $\rho$  \endsubhead 
 Let us introduce the following function
 $$\rho: {\Bbb R} \longrightarrow [0, \pi], \quad \rho(x)=\min_{k \in {\Bbb Z}} |x- 2\pi k|.$$
 In words: $\rho(x)$ is the distance from $x$ to the closest integer multiple of $2\pi$.
 Clearly,
 $$\rho(-x)=\rho(x) \quad \text{and} \quad \rho(x + y) \ \leq \ \rho(x) +\rho(y)$$
 for all $x, y \in {\Bbb R}$. 
We will use that 
$$1-{1 \over 2} \rho^2(x)\ \leq \ \cos x  \ \leq \ 1-{1 \over 5}\rho^2(x). \tag7.2.1$$

\subhead (7.3) The absolute value of $F(s,t)$ \endsubhead
Let 
$$\alpha_{jk} =2 \zeta_{jk} \left(1 + \zeta_{jk} \right) \quad \text{for all} \quad j, k.$$
Then 
$$2\delta^2 \tau^2  \ \leq \ \alpha_{jk} \ \leq \ 4 \tau^2 \quad 
\text{for all} \quad j, k \tag7.3.1$$

Let us define functions
$$f_{jk}(x)={1 \over \sqrt{1+\alpha_{jk} -\alpha_{jk} \cos x}} \quad \text{for} \quad x \in {\Bbb R}. $$
Then we can write
$$ |F(s,t)|=
\prod \Sb 1 \leq j \leq m \\ 1\leq k \leq n \endSb f_{jk}\left(s_j +t_k\right).\tag7.3.2$$
We observe that 
$$f_{jk}(x) =1 \quad \text{provided} \quad \rho(x)=0 \quad \text{and that} \quad 
f_{jk}(-x)=f_{jk}(x).$$
From (7.2.1) and (7.3.1) we conclude that for any $\epsilon >0$ 
we have 
$$\aligned f_{jk}(x) \ \leq \ &\exp\left\{ -\gamma(\delta, \epsilon)\right\} f_{jk}(y) \\
&\text{provided} \quad \rho(x) \ \geq \ {2 \epsilon \over \tau} \quad \text{and} \quad 
\rho(y) \ \leq \ {\epsilon \over \tau}, \endaligned \tag7.3.3$$
where $\gamma(\delta, \epsilon)>0$ is a constant.

Finally, we observe that 
$${d^2 \over dx^2} \ln f_{jk}(x)=
-{\alpha_{jk}(1+\alpha_{jk}) \cos x -\alpha_{jk}^2 \over 2(1+\alpha_{jk}-\alpha_{jk} \cos x)^2}.$$
It follows from (7.2.1) and (7.3.1) that for all $j, k$
$${d^2 \over d x^2} \ln f_{jk}(x) \ \leq \ - {2 \over 5} \tau^2 \delta^2 \quad \text{provided} 
\quad |x| \leq {\delta^2 \over 5 \tau}. \tag7.3.4$$
In particular, the function $\ln f_{jk}(x)$ is strictly log-concave on the interval 
$|x| \leq \delta^2/5 \tau$.
\bigskip
In what follows, we fix a particular parameter $\epsilon >0$. All implied constants in the ``$O$'' and ``$\Omega$'' notation may depend only 
on the parameters $\delta$ and $\epsilon$. We say that $m$ and $n$ are {\it sufficiently large} if 
$m +n \geq \gamma(\delta, \epsilon)$ for some constant $\gamma(\delta, \epsilon)>0$.
Recall that $m$ and $n$ are of the same order, $m \geq \delta n$ an $n \geq \delta m$.

Our first goal is to show that for any fixed $\epsilon >0$ only the points $(s, t) \in \Pi_0$ 
for which the inequality $\rho\left(s_j + t_k \right) \leq \epsilon/\tau$ holds for the the overwhelming majority of pairs of indices $(j, k)$ contribute significantly to the integral of $|F(s, t)|$ on $\Pi_0$.
Recall that $\tau_n =0$ on $\Pi_0$.

\proclaim{(7.4) Lemma} For $\epsilon >0$ and a point $(s, t) \in \Pi_0$ let us define 
the following two sets:
\smallskip
Let $J=J(s, t; \epsilon) \subset \{1, \ldots, m\}$ be the set of all indices $j$ 
such that 
$$\rho\left(s_j + t_k \right) \ \leq \ \epsilon/\tau$$ for more than $(n-1)/2$ distinct indices $k=1, \ldots, n-1$

and

Let $K=K(s, t; \epsilon) \subset \{1, \ldots, n-1\}$ be the set of all indices $k$ for 
such that 
$$\rho\left(s_j +t_k \right) \ \leq \ \epsilon/\tau$$ for more than $m/2$ distinct indices $j=1, \ldots, m$.
\smallskip
Let $\overline{J}=\left\{1, \ldots, m \right\} \setminus J$ and let 
$\overline{K} =\left\{1, \ldots, n-1\right\} \setminus K$.

Then 
\roster
\item
$$\split &|F(s, t)| \  \leq \  \exp\left\{ -\gamma(\delta, \epsilon) n |\overline{J}|  \right\} \quad 
\text{and} \\ &|F(s, t)| \ \leq \ \exp\left\{-\gamma(\delta, \epsilon) m |\overline{K}| \right\} \endsplit$$
for some constant $\gamma(\delta, \epsilon)>0$.
\item
$$\split &\rho\left(s_{j_1} - s_{j_2} \right)\ \leq \ 2\epsilon/\tau \quad \text{for all} \quad j_1, j_2 \in J 
\quad \text{and}\\  &\rho\left(t_{k_1} - t_{k_2} \right)\ \leq \ 2\epsilon/\tau \quad 
\text{for all} \quad k_1, k_2 \in K. \endsplit$$
\item If $|J| > m/2$ or $|K| > (n-1)/2$ then 
$$\rho\left(s_j + t_k \right) \ \leq \ 3\epsilon/\tau \quad \text{for all} \quad j \in J \quad \text{and all} \quad k \in K.$$
\endroster
\endproclaim
\demo{Proof}
For every $j \in \overline{J}$ there are at least $(n-1)/2$ distinct $k$ 
and for every $k \in \overline{K}$ there are at least $m/2$ distinct $j$
such that 
$$\rho\left(s_j + t_k \right) > \epsilon/\tau$$
and hence 
$$f_{jk}\left(s_j + t_k\right) \ \leq \ \exp\left\{-\Omega(1) \right\}$$
by (7.3.3).
Part (1) follows from (7.3.2).

For every $j_1, j_2 \in J$ there is at least one common index $k$ such that 
$$\rho\left( s_{j_1} + t_k \right) \leq \epsilon/\tau 
\quad \text{and} \quad  \rho\left(s_{j_2} + t_k \right) \leq \epsilon/\tau.$$
Then 
$$\rho\left(s_{j_1} -s_{j_2} \right) = \rho\left(s_{j_1} +t_k -s_{j_2} - t_k \right) 
\ \leq \ \rho\left(s_{j_1} + t_k \right) + \rho\left(s_{j_2} + t_k \right) \ \leq \ 2 \epsilon/\tau.$$
The second inequality of Part (2) follows similarly.

If $|K| > (n-1)/2$ then for every $j \in J$ there is a $k_j \in K$ such that 
$$\rho\left(s_j + t_{k_j} \right) \leq \epsilon/\tau.$$
Then, by Part (2), for every $k \in K$ we have 
$$\rho\left(s_j +t_k \right) =\rho\left(s_j + t_{k_j} -t_{k_j} + t_k \right) 
\ \leq \ \rho\left(s_j +t_{k_j} \right) + \rho\left(t_k -t_{k_j}\right)\ \leq \ 3 \epsilon/\tau.$$
The case of $|J| > m/2$ is handled similarly.
{\hfill \hfill \hfill} \qed
\enddemo

Using estimates of Theorem 6.1, it is not hard to deduce from Part (1) of Lemma 7.4 that for 
any fixed $\epsilon >0$ only points $(s, t) \in \Pi_0$ with 
$\overline{J}=O(\ln m)$ and $\overline{K}=O(\ln n)$ may contribute essentially to the integral 
of $|F(s, t)|$. It follows then by Part (3) of Lemma 7.4 that for such points we have 
$\rho\left(s_j +t_k \right) \leq 3 \epsilon/\tau$ for all $j \in J$ and all $k \in K$. 
 Our next goal is to show that only those points $(s, t) \in \Pi_0$ contribute 
substantially to the integral for which $\rho\left(s_j +t_k \right) \leq \epsilon/\tau$ {\it for all}
pairs $(j, k)$.

\proclaim{(7.5) Proposition} For $\epsilon >0$ let us define a set 
$X(\epsilon) \subset \Pi_0$ by 
$$\split X(\epsilon)=\Bigl\{(s, t) \in \Pi_0: \quad &\rho\left(s_j +t_k \right) \leq \epsilon/\tau
\\  &\text{for all} \quad j=1, \ldots, m \quad \text{and} \quad k=1, \ldots, n-1 \Bigr\}. \endsplit$$
Then
$$\int_{\Pi_0 \setminus X(\epsilon)} |F(s,t)| \ ds dt \ \leq \ \exp\bigl\{-\gamma(\delta, \epsilon)
(m+n)\bigr\}
 \int_{\Pi_0} |F(s,t)| \ ds dt$$
 for some constant $\gamma(\delta, \epsilon)>0$ and all sufficiently large $m+n$.
\endproclaim
\demo{Proof} 
For subsets $A \subset \{1, \ldots, m\}$ and $B \subset \{1, \ldots, n-1\}$
let us define 
a set $P_{A, B} \subset \Pi_0$ (we call it a {\it piece})  by 
$$P_{A, B} =\Bigl\{(s, t) \in \Pi_0: \quad \rho\left(s_j+ t_k \right) \leq \epsilon/40\tau \quad \text{for all}
 \quad 
j \in A \quad \text{and all} \quad k \in B\Bigr\}.$$
Let 
$$V= \bigcup_{A, B} P_{A, B},$$
where the union is taken over all subsets $A$ and $B$ such that
$$|\overline{A}| \leq \ln^2 m \quad \text{and} \quad \overline{B} \leq \ln^2 n,$$
where
 $$\overline{A}=\{1, \ldots, m\} \setminus A \quad \text{and} \quad \overline{B}=\{1, \ldots, n-1\}\setminus B.$$
We claim that the integral over $\Pi_0 \setminus V$ is asymptotically negligible.
Indeed, for sufficiently large $m+n$ and
for all $(s, t) \in \Pi \setminus V$, by Part (3) of Lemma 7.4, we must have 
$$\overline{J}(s, t; \epsilon/120) \ \geq \ \ln^2 m \quad \text{or} \quad 
\overline{K}(s, t; \epsilon/120) \ \geq \ \ln^2 n.$$ 
In either case, by Part (1) of Lemma 7.4, we must have 
$$|F(s,t)| \ \leq \ \exp\left\{ -\Omega\bigl(n \ln^2 n \bigr)\right\}$$
for all sufficiently large $m+n$.
On the other hand, by Parts (1), (2) and (4) of Theorem 6.1, we conclude that 
$$\int_{\Pi_0} |F(s, t)| \ ds dt \ \geq \ \exp\left\{ - O\bigl(n \ln n \bigr)\right\}$$
(we use that $\tau$ is bounded by a polynomial in $m+n$).
This proves that 
$$\int_{\Pi_0 \setminus V} |F(s,t)| \ ds dt \ \leq \ \exp\left\{  -\Omega\left(n \ln^2 n\right) 
\right\} \int_{\Pi_0} | F(s,t)| \ ds dt \tag7.5.1$$
provided $m+n$ is sufficiently large.

Next, we prove that the integral over $V\setminus X(\epsilon)$ is asymptotically negligible.

As in the proof of Part (2) of Lemma 7.4, we conclude that for every piece $P_{A, B}$ and 
for every $(s,t) \in  P_{A,B}$ we have
$$\aligned &\rho\left(s_{j_1}-s_{j_2} \right) \ \leq \  \epsilon/20\tau  \quad \text{for all} \quad j_1, j_2 \in A \qquad
\text{and} \\
&\rho\left(t_{k_1} -t_{k_2} \right) \ \leq \  \epsilon/20\tau  \quad \text{for all} \quad k_1, k_2 \in B. \endaligned
\tag7.5.2$$
 
 Let us choose a point $(s, t) \in P_{A, B} \setminus X(\epsilon)$.
  Hence we have 
 $\rho\left(s_{j_0} +t_{k_0} \right) > \epsilon/\tau$ for some $j_0$ and $k_0$. Let us pick 
 any $j_1 \in A$ 
 and $k_1 \in B$. Then 
 $$\split \rho\left(s_{j_0} +t_{k_0}\right)=&\rho\left(s_{j_0} + t_{k_0} + s_{j_1} + t_{k_1} -s_{j_1} -t_{k_1} \right) \\
 \ \leq \ &\rho\left(s_{j_0} + t_{k_1} \right) + \rho\left(s_{j_1} + t_{k_0} \right) + 
 \rho\left(s_{j_1} +t_{k_1} \right). \endsplit$$
 Since $\rho\left(s_{j_1} + t_{k_1}\right) \leq \epsilon/40 \tau$, we must have 
 either 
 $$\rho\left(s_{j_0} +t_{k_1} \right)\  >\  39\epsilon/80 \tau,$$ 
 in which case necessarily $j_0 \in \overline{A}$, or 
 $$\rho\left(s_{j_1}+t_{k_0}\right) \ > \ 39 \epsilon/80\tau,$$
 in which case necessarily $k_0 \in \overline{B}$. 
 
 In the first case (7.5.2) implies that 
 $$\rho\left(s_{j_0} +t_k \right)\ >  \ 35\epsilon/80\tau = 7\epsilon/16\tau \quad \text{for all}\quad  k \in B$$
  and in the second case (7.5.2) 
 implies that 
 $$\rho\left(s_j + t_{k_0} \right) > 35 \epsilon/80\tau =7\epsilon/16\tau  \quad \text{for all} \quad j \in A.$$
 
 For $j \in \overline{A}$ we define
 $$Q_{A, B; j}=\Bigl\{(s,t) \in P_{A, B}: \quad \rho\left(s_j +t_k \right) >  7\epsilon/16\tau 
 \quad \text{for all} \quad k \in B\Bigr\} \tag7.5.3$$
 and for $k \in \overline{B}$ we define
 $$R_{A, B; k}=\Bigl\{ (s, t) \in P_{A,B}: \quad \rho\left(s_j +t_k \right) > 7\epsilon/16\tau 
 \quad \text{for all} \quad j \in A \Bigr\}. \tag7.5.4$$
 Then
 $$P_{A,B}\setminus X(\epsilon) \quad \subset \quad \left(\bigcup_{j \in \overline{A}} Q_{A, B; j} \quad \bigcup \quad
 \bigcup_{k \in \overline{B}} R_{A, B; k}\right). \tag7.5.5$$ 
 Let us compare the integrals
  $$\int_{Q_{A, B; j}} |F(s, t)| \ ds dt \quad \text{and} \quad \int_{P_{A, B}} |F(s, t)| \ ds dt.$$
 Given a point $(s, t) \in P_{A, B}$ we obtain another 
point in $P_{A, B}$ if we arbitrarily choose coordinates $s_j \in [-\pi, \pi]$ for $j \in \overline{A}$
 and $t_k \in [-\pi, \pi]$ for $k \in \overline{B}$. 
 Let us pick a particular non-empty set $Q_{A, B; j_0}$ for some $j_0 \in \overline{A}$.
 We obtain a {\it fiber} $E =E_{j_0} \subset P_{A,B}$ if 
we let the coordinate $s_{j_0}$ vary 
arbitrarily between $-\pi$ and $\pi$ while fixing all other coordinates of some point $(s, t) \in P_{A, B}$.
Geometrically, each fiber $E$ is an interval of length $2\pi$. We construct a set $I \subset E$ 
as follows: we choose an arbitrary coordinate $k_1 \in B$ and let $s_{j_0}$ vary in such a way 
that $\rho\left(s_{j_0} + t_{k_1} \right) \leq  \epsilon/20\tau$. Geometrically, $I$ is an interval of 
length $ \epsilon/10\tau$ or a union of two non-overlapping intervals of the total length 
$ \epsilon/10 \tau$.
Moreover, by (7.5.2), we have 
$$\rho\left(s_{j_0} + t_k \right) \ \leq \  \epsilon/10\tau \quad \text{for all} \quad k \in B \quad 
\text{and all} \quad (s, t) \in I. \tag7.5.6$$

As we vary $s_{j_0}$ without changing other coordinates, in the product 
(7.3.2) only the functions $f_{j_0 k}$ change.
 Comparing (7.5.6) and (7.5.3) and using (7.3.2) and (7.3.3),
  we conclude that
  $$\split &|F(s, t)| \ \leq \ \exp\left\{-\Omega(n)\right\} \left|F\left(\tilde{s}, \tilde{t}\right)\right|
  \\ &\quad \text{for all} \quad (s, t) \in Q_{A, B; j_0} \cap E  \quad \text{and all} \quad
\left(\tilde{s}, \tilde{t} \right) \in I. \endsplit$$
Therefore, 
$$\int_{E \cap Q_{A,B; j_0} } |F(s,t)| \ ds_{j_0} \ \leq \  
\exp\left\{-\Omega( n) \right\}  \int_E |F(s,t)| \ ds_{j_0}$$
 provided $m+n$ is large enough
 (again, we use that $\tau$ is bounded by a polynomial in $m+n$).
 Integrating over all fibers $E \subset P_{A, B}$, we prove that  
$$\int_{Q_{A,B; j} } |F(s,t)| \ ds dt \ \leq \  \exp\left\{ -\Omega(n) \right\} \int_{P_{A,B}}|F(s,t)| \ ds dt$$
  provided $m+n$ is large enough. Similarly, we prove that for sets $R_{A,B; k}$ defined by (7.5.4)
  we have 
$$\int_{R_{A,B; k} } |F(s,t)| \ ds dt \ \leq \  \exp\left\{-\Omega(n)\right\} \int_{P_{A,B}}|F(s,t)| \ ds dt$$
 provided $m+n$ is large enough.  
Since $|\overline{A}| \leq \ln^2 m$ and $\overline{B}| \leq \ln^2 n$, from (7.5.5) we deduce that
  $$\int_{P_{A, B} \setminus X(\epsilon)} 
  |F(s,t)| \ ds dt \ \leq \  \exp\left\{ -\Omega(n)\right\} \int_{P_{A,B}}|F(s,t)| \ ds dt$$
 provided $m+n$ is large enough. Finally, since the number of pieces 
 $P_{A,B}$ does not exceed $\exp\left\{O\left( \ln^3 n \right)\right\}$,
 the proof follows by (7.5.1).
 {\hfill \hfill \hfill} \qed
\enddemo

Our next goal is to show that the integral over $\Pi_0 \setminus I$ is negligible, where 
$$I=\Bigl\{(s, t) \in \Pi_0: \quad |s_j|, |t_k| \ \leq \ \epsilon/\tau \quad \text{for all} \quad j, k \Bigr\}.$$
We accomplish this in the next two lemmas.

\proclaim{(7.6) Lemma} For $\epsilon>0 $ let us define a set $Y(\epsilon) \subset \Pi_0$ 
by 
$$\split Y(\epsilon)=\Bigl\{(s, t) \in \Pi_0: \quad & |s_j +t_k| \leq \epsilon/\tau  \\& \text{for all} \quad
j=1, \ldots, m \quad \text{and all} \quad k=1, \ldots, n-1\Bigr\}. \endsplit$$
Then 
$$\int_{\Pi_0 \setminus Y(\epsilon)} |F(s, t)| \ ds dt \ \leq \ \exp\bigl\{-\gamma(\delta, \epsilon)(m+n) \bigr\}
\int_{\Pi_0} |F(s, t)| \ ds dt,$$
for some constant $\gamma(\epsilon, \delta) > 0$ and all sufficiently large $m+n$.
\endproclaim 
\demo{Proof} Without loss of generality, we assume that $\epsilon <\delta/2$, so $\epsilon/\tau <1/2$.

Let $X(\epsilon)$ be the set of Proposition 7.5 and let us define 
$$\split &Z(\epsilon) =\Bigl\{(s, t) \in \Pi_0: \quad s_j, t_k \in [-\pi, \ -\pi+\epsilon/\tau] \cup [\pi-\epsilon/\tau, \ \pi] \\
&\qquad \qquad \qquad \text{for all} \quad j=1, \ldots, m \quad \text{and all} \quad 
k=1, \ldots, n-1\Bigr\}. \endsplit$$
We claim that 
$$X(\epsilon) \ \subset \ Y(\epsilon) \ \cup \ Z(2\epsilon). \tag7.6.1$$
We note that if $ \rho(x) \ \leq \ \epsilon/\tau$ for some $-2\pi \ \leq x \ \leq \ 2\pi$ then either $|x| \leq \epsilon/\tau$ or $x \geq 2\pi-\epsilon/\tau$ or $x \leq -2\pi +\epsilon/\tau$. To prove (7.6.1),
let us pick an arbitrary $(s, t) \in X(\epsilon)$. Suppose that
$$-\pi + 2\epsilon/\tau\ < \ s_{j_0} \ < \ \pi -2 \epsilon/\tau \quad \text{for some} \quad j_0. \tag7.6.2$$
Since $- \pi \leq t_k \leq \pi$ for all $k$, we have 
$$-2 \pi + 2\epsilon/\tau \ < \ s_{j_0} + t_k  \ < \ 2\pi-2 \epsilon/\tau
\quad \text{for} \quad k=1, \ldots, n.$$
Since $\rho\left(s_{j_0} +t_k \right) \leq \epsilon/\tau$, we must have 
$$\left| s_{j_0} +t_k \right| \ \leq \ \epsilon/\tau \quad \text{for} \quad k=1, \ldots, n-1$$
and, therefore,
$$-\pi +\epsilon/\tau \ < \ t_k \ < \ \pi-\epsilon/\tau \quad \text{for} \quad k=1, \ldots, n-1.$$
Since $-\pi \leq s_j \leq \pi$ we conclude that 
$$\split &-2\pi +\epsilon/\tau \ < \ s_j + t_k \ < \ 2\pi -\epsilon/\tau \\ &\qquad \qquad \text{for all}
\quad j=1, \ldots, m \quad \text{and} \quad k=1, \ldots, n-1. \endsplit$$
Since $\rho\left(s_j +t_k \right) \leq \epsilon/\tau$ we conclude that 
$$\split &\left| s_j +t_k\right| \ \leq \ \epsilon/\tau \\ &\qquad \qquad  \text{for all} \quad
j=1, \ldots, m \quad \text{and} \quad k=1, \ldots, n-1 \endsplit$$
and hence $(s, t) \in Y(\epsilon)$.

Similarly, we prove that if 
$$-\pi + 2\epsilon/\tau \ < \ t_{k_0} \ < \ 2 \pi -2 \epsilon/\tau \quad \text{for some} \quad 
k_0 \tag7.6.3$$ then 
$(s, t)   \in Y(\epsilon)$. If both (7.6.2) and (7.6.3) are violated, then $(s, t) \in Z(2\epsilon)$ and so we obtain 
(7.6.1).

Next, we show that the integral over $Z(2\epsilon)$ is asymptotically negligible.
The set $Z(2\epsilon)$ is a union of $2^{m+n-1}$ pairwise disjoint
{\it corners}, where each corner is determined by a choice of the interval $[-\pi, -\pi +2\epsilon/\tau]$
or $[\pi-2\epsilon/\tau, \pi]$ for each coordinate $s_j$ and $t_k$.
The transformation
$$\split &s_j \longmapsto \cases s_j + \pi &\text{if\ } s_j \in [-\pi, \ -\pi +2\epsilon/\tau] 
\\ s_j -\pi &\text{if\ }
s_j \in [\pi-2 \epsilon/\tau, \ \pi] \endcases \qquad \text{for} \quad j=1, \ldots, m \\
 &\qquad \qquad \text{and} \\
&t_k \longmapsto \cases t_k+ \pi &\text{if\ } t_k \in [-\pi, \ -\pi +2\epsilon/\tau] \\ t_k -\pi &\text{if\ }
t_k \in [\pi-2 \epsilon/\tau, \ \pi] \endcases \qquad \text{for} \quad k=1, \ldots, n-1\endsplit$$
is measure-preserving and maps $Z(2\epsilon)$ onto the cube 
$$I=\Bigl\{\left(s, t\right): \quad \left|s_j\right|, \ \left|t_k\right| 
\ \leq \ 2 \epsilon/\tau \quad \text{for all} \quad j, k \Bigr\}. $$
In the product (7.3.2), it does not change the value of $f_{jk}$ except when $k=n$ (recall that 
$t_n=0$ on $\Pi_0$).
Since $2\epsilon/\tau < 1$, by (7.3.3) the transformation increases the value of each function $f_{jn}$ 
by at least a
factor of $\gamma(\delta) >1$. Therefore,
$$\int_{Z(2\epsilon)} |F(s, t)| \ ds dt \ \leq \ \exp\left\{-\Omega(m) \right\} \int_I |F(s, t)| \ ds dt$$
and the proof follows by (7.6.1) and Proposition 7.5.
{\hfill \hfill \hfill} \qed
\enddemo

\proclaim{(7.7) Lemma} For $\epsilon >0$ let us define the cube
$$I(\epsilon)=\Bigl\{(s,t) \in \Pi_0: \quad |s_j|, |t_k| \leq \epsilon/\tau \quad \text{for all} \quad
j, k \Bigr\}.$$
Then
$$\int_{\Pi_0 \setminus I(\epsilon)} |F(s,t)| \ ds dt \ \leq \ \exp\bigl\{-\gamma(\delta, \epsilon) (m+n)
\bigr\} \int_{\Pi_0} |F(s,t)| \ ds dt$$
for some $\gamma(\delta, \epsilon)>0$ and $m+n$ large enough.
\endproclaim 
\demo{Proof} Without loss of generality, we assume that $\epsilon < \delta$, so $\epsilon/\tau <1$.

Let $Y(\epsilon/20)$ be the set of Lemma 7.6, so the integral of $|F(s, t)|$ over \break
$\Pi \setminus Y(\epsilon/20)$ is asymptotically negligible.

Let us choose a point $(s, t) \in Y(\epsilon/20)$. We have 
$$l \epsilon/20\tau \ \leq \ s_1\  \leq \ (l+1) \epsilon/20\tau \qquad \text{for some integer} \quad l.$$
Since $|s_1 +t_k| \leq \epsilon/20$, we obtain
$$(-l-2)\epsilon/20\tau \ \leq \ t_k \ \leq (-l+1) \epsilon/20\tau \quad \text{for} \quad k=1, \ldots, n-1$$
and then similarly
$$(l-2) \epsilon/20\tau  \ \leq \ s_j \ \leq \ (l+3) \epsilon/20\tau \quad \text{for}  \quad j=1,\ldots, m.$$
Let us denote 
$$w=\left(\underbrace{\epsilon/20\tau, \ldots, \epsilon/20\tau}_{\text{$m$ times}};  \ 
\underbrace{-\epsilon/20\tau, \ldots, -\epsilon/20\tau}_{\text{$n-1$ times}}, 0 \right).$$
Hence we conclude that 
$$Y(\epsilon/20) \ \subset \ \bigcup\Sb |l| \ \leq 1+ 20 \pi \tau/\epsilon \\ l \in {\Bbb Z} \endSb
   I(3\epsilon/20) +lw. \tag7.7.1$$
Since 
$\tau$ is bounded by a polynomial in $m$ and $n$, the number 
of translates of the cube $I(3 \epsilon/20)$ in the right hand side of (7.7.1) is $(m+n)^{O(1)}$.

The translation
$$(s, t) \longmapsto (s, t) +l w$$
does not change the value of the functions $f_{jk}(s_j +t_k)$ in (7.3.2), unless $k=n$
(recall that $t_n=0$ on $\Pi_0$).

 For $(s, t) \in I(3\epsilon/20)$ we have $|s_j| \leq 3\epsilon/20\tau$ for all $j$.
 For $(s, t) \in I(3\epsilon/20) + lw$ with $|l| \geq 10$, we have $|s_j| \geq 7\epsilon/20\tau$ for all $j$.
 Since $\epsilon/\tau <1$, for all $l$ in the union of (7.7.1) such that $|l| \geq 10$ and all 
 $(s, t) \in I(3 \epsilon/20) +lw$ we have 
 $\rho(s_j) \geq 6 \epsilon/20\tau$ for all $j=1, \ldots, m$. 
 
 Using (7.3.2) and (7.3.3) we conclude that 
 $$\split  |F(s, t)| \ \leq \ \exp\left\{ - \Omega(m) \right\} &|F(\tilde{s}, \tilde{t})| 
 \quad \text{for all} \quad (\tilde{s}, \tilde{t}) \in I(3 \epsilon/20) \\
 & \text{and for all} \quad (s, t) \in I(3 \epsilon/20) +l w \quad \text{with} \quad |l| \geq 10. \endsplit$$
Since the number of translates in (7.7.1) is bounded by a polynomial in $(m+n)$ and since 
$$I(3 \epsilon/20) +l w \ \subset \ I(\epsilon) \quad \text{provided} \quad |l| \leq 10,$$
the proof follows by Lemma 7.6.
 {\hfill \hfill \hfill} \qed
 \enddemo

To finish the proof of Theorem 7.1 we need a concentration inequality for strictly log-concave
probability measures.  

\proclaim{(7.8) Theorem} Let $V$ be Euclidean space with the norm 
$\| \cdot\|$,  let  $B \subset V$ be a convex body,  let us consider a probability 
measure supported on $B$ with the density $e^{-U}$, where $u: B \longrightarrow {\Bbb R}$ 
is a function satisfying 
$$U(x) +U(y) -2U \left({x+y \over 2} \right) \ \geq \ c \|x-y\|^2 \quad \text{for all} \quad x, y \in B$$
and some constant $c>0$. For a point $x \in V$ and a closed subset $A \subset V$ 
we define the distance
$$\dist(x, A)=\min_{y \in A} \|x-y\|.$$
Let $A \subset B$ be a closed set such that $\PP(A) \geq 1/2$. Then,
for any $r \geq 0$ we have 
$$\PP\Bigl\{x \in B: \quad \dist(x, A) \geq r \Bigr\} \ \leq \ 2 e^{-cr^2}.$$
\endproclaim
\demo{Proof} See, for example,
Section 2.2 of \cite{Le01} or Theorem 8.1 and its proof in \cite{B97a},  which, although 
stated for the Gaussian measure is adapted in a straightforward way to our situation.
\enddemo
{\hfill \hfill \hfill} \qed

Here is how we apply Theorem 7.8.
\proclaim{(7.9) Lemma} Let us choose $0 < \epsilon \leq \delta^2/10$.
In the space ${\Bbb R}^{m+n}$ let us consider the hyperplane 
$$H=\left\{\left(s_1, \ldots, s_m; t_1, \ldots, t_n \right): \quad \sum_{j=1}^m s_j 
=\sum_{k=1}^n t_k\right\}.$$
Let $B \subset H$ be a convex body centrally symmetric about the origin:
$(s, t) \in B$ if and only if $(-s, -t) \in B$, and such that for all $(s, t) \in B$ we have
$$\split &|s_j| \leq \epsilon/\tau \quad \text{for} \quad j=1, \ldots, m \\
&|t_k| \leq \epsilon/\tau \quad \text{for} \quad k=1, \ldots, n. \endsplit$$

Let us consider the probability measure on $B$ with the density proportional to $|F(s,t)|$.
Then, for any $\kappa>0$ we have 
$$\PP\Bigl\{ (s,t) \in B: \quad |s_j|, |t_k| \ \leq \ {\ln(m+n) \over 2\tau \sqrt{m+n}}
 \quad \text{for all} \quad j, k \Bigr\} 
\ \geq \ 1 - (m+n)^{-\kappa},$$
provided $m+n \geq \gamma(\delta, \kappa)$ for some constant $\gamma(\delta, \kappa)>0$.
\endproclaim 
\demo{Proof} Let $f_{jk}$ be the functions defined in Section 7.3 and let 
$$u_{jk} = -\ln f_{jk}.$$
We define 
 $$U(s, t) =a+\sum \Sb 1 \leq j \leq m \\ 1 \leq k \leq n \endSb u_{jk} \left(s_j + t_k \right) \
\quad \text{for} \quad (s, t) \in B,$$
where $a$ is a constant chosen in such a way that 
$$e^{-U}=e^{-a} |F(s, t)|$$
is a probability density on $B$. It follows by (7.3.4) that 
$$\split &u_{jk}(x) +u_{jk}(y) -2 u_{jk}\left({x+y \over 2} \right) \ \geq \ \Omega\left( \tau^2 (x-y)^2 \right) \\
&\qquad \text{provided} \quad |x|, |y| \ \leq \ 2\epsilon/\tau. \endsplit \tag7.9.1 $$ 
Let us consider the map 
$$M: \left(s_1, \ldots, s_m; t_1, \ldots, t_n \right) \longmapsto \left( \ldots s_j +t_k \ldots \right)$$
as a map $M: H \longmapsto {\Bbb R}^{mn}$. 
From Lemma 3.6
$$\|M x\|^2  \ \geq \ \min\{m, n\} \|x\|^2  \quad \text{for all} \quad x \in H,$$
where $\| \cdot \|$ is the Euclidean norm in the corresponding space.
It follows then by (7.9.1) that 
 $$U(x) +U(y) -2 U\left({x+y \over 2} \right) \ \geq \ \Omega\left( \tau^2 n \|x-y\|^2\right) 
 \quad \text{for all} \quad x, y \in B.$$
 Now we apply Theorem 7.8 with
 $$c=  \Omega\left( \tau^2 n \right)$$
 to the probability density $e^{-U}$ on $B$.

For $j=1, \ldots, m$, let $S_j^+$ be the set consisting of the points $(s, t) \in B$ with $s_j \geq 0$,
let $S_j^-$ be the set consisting of the points $(s, t) \in B$ with $s_j \leq 0$, let $T_k^+$ be the set 
consisting of the points $(s, t) \in B$ with $t_k \geq 0$ and let $T_k^-$ be the set consisting of the points 
$(s, t) \in B$ with $t_k \leq 0$. Since both $B$ and the probability measure are invariant under the symmetry 
$$(s, t) \longmapsto (-s, -t),$$
we have 
$$\PP\left(S_j^+\right)=\PP\left(S_j^-\right)= \PP\left(T_k^+\right)= \PP\left(T_k^-\right) = {1 \over 2}.$$
We note that 
$$\split &\text{if} \quad \left|s_j\right| \geq r \quad \text{then} \quad \dist\bigl((s, t),\ S_j^+\bigr), \ 
\dist\bigl((s, t),\ S_j^- \bigr) \geq r \qquad \text{and} \\
&\text{if} \quad \left|t_k\right| \geq r \quad \text{then} \quad \dist\bigl((s, t),\ T_k^+\bigr), \ 
\dist\bigl((s, t),\ T_k^- \bigr) \geq r. \endsplit$$
Applying Theorem 7.8 with  
$$r={\ln (m+n) \over 2\tau \sqrt{m+n}}$$
We conclude that for all $j$ and $k$ 
$$\split &\PP\left\{ (s, t) \in B: \quad |s_j| > {\ln (m+n) \over 2 \tau \sqrt{m+n}} \right\} \ \leq \ 
\exp\left\{ - \Omega\left(\ln^2 m \right) \right\} \quad \text{and} \\
 &\PP\left\{ (s, t) \in B: \quad |t_k| > {\ln (m+n) \over 2 \tau \sqrt{m+n}} \right\} \ \leq \ 
\exp\left\{ - \Omega\left(\ln^2 n \right) \right\} \endsplit$$
and the proof follows.
{\hfill \hfill \hfill} \qed
\enddemo

Now we are ready to prove Theorem 7.1.

\subhead (7.10) Proof of Theorem 7.1 \endsubhead Let us choose an $0< \epsilon \leq \delta^2/10$ 
as in Lemma  7.9 and let $H \subset {\Bbb R}^{m+n}$ be the hyperplane 
defined in Lemma 7.9. We identify ${\Bbb R}^{m+n-1}$ with the hyperplane $\tau_n=0$ in 
${\Bbb R}^{m+n}$. 
We consider a linear transformation
$T:\ H \longrightarrow {\Bbb R}^{m+n-1}$, 
$$ \left(s_1, \ldots, s_m;\ t_1, \ldots, t_n \right) \longmapsto 
\left(s_1+t_n, \ldots, s_m +t_n;\ t_1-t_n, \ldots, t_{n-1}-t_n, 0 \right).$$
The inverse linear transformation 
$$T^{-1}: \quad \left(s_1', \ldots, s_m';\ t_1', \ldots, t_{n-1}', 0 \right) \longmapsto 
\left(s_1, \ldots, s_m;\ t_1, \ldots, t_n \right)$$
is computed as follows:
$$t_n={s_1' + \ldots + s_m' -t_1' - \ldots - t_{n-1}' \over m+n}, \quad s_j=s_j'-t_n, \quad t_k=t_k'+t_n.$$
Let us consider the cube $I=I(\epsilon/2) \subset {\Bbb R}^{m+n-1}$ defined by the 
inequalities
$$|s_j|,\  |t_k| \ \leq \ \epsilon/2 \tau \quad \text{for} \quad j=1, \ldots, m \quad \text{and} \quad
k=1, \ldots, n-1.$$
By Lemma 7.7 we have 
$$\int_{\Pi_0 \setminus I} |F(s, t)| \ ds dt \ \leq \ \exp\left\{-\Omega(m+n) \right\} \int_{\Pi_0} | F(s, t)| \ ds dt
\tag7.10.1$$
for all sufficiently large $m$ and $n$. 

Let $B=T^{-1}(I) \subset H$. Then $B$ is centrally symmetric and convex,  and 
for all $(s, t) \in B$ we have 
$|s_j|, |t_k| \leq \epsilon/\tau$ for all $j$ and $k$. Let 
$$\split A=\Bigl\{(s, t) \in B: \quad &|s_j| \leq {\ln(m+n) \over 2\tau\sqrt{m+n}} \quad \text{for} \quad j=1,
\ldots, m
\qquad \text{and} \\
&|t_k|  \leq {\ln(m+n) \over 2\tau\sqrt{m+n}} \quad \text{for} \quad k=1, \ldots, n \Bigr\}. \endsplit$$
By Lemma 7.9,
$$\int_{B \setminus A} |F(s, t)| \ ds dt \ \leq \ (m+n)^{-\kappa} \int_B |F(s, t)| \ ds dt$$
provided $m+n \geq \gamma(\delta, \kappa)$ for some $\gamma(\delta, \kappa) >0$. 
Now, the push-forward of the probability measure on $B$ with the density 
proportional to $|F(s, t)|$ under the map $T$ is the probability measure on $I$ with the density proportional 
to $|F(s, t)|$. Moreover, the image $T(A)$ lies in the cube $\UU$ defined by the inequalities
$$|s_j|\ |t_k|  \ \leq \ {\ln (m+n) \over \tau \sqrt{m+n}} \quad \text{for} \quad j=1, \ldots, m
\quad \text{and} \quad k=1, \ldots, n-1.$$
Therefore, 
$$\int_{I \setminus \UU} |F(s, t) \ ds dt \ \leq \ (m+n)^{-\kappa} \int_I |F(s, t)| \ ds dt. \tag7.10.2$$
provided $m+n \geq \gamma(\delta, \kappa)$ for some $\gamma(\delta, \kappa)>0$.
The proof now follows by (7.10.1) and  (7.10.2).
{\hfill \hfill \hfill} \qed

\head 8. Proof of Theorem 1.3  \endhead

First, we prove Theorem 1.3 assuming, additionally, that 
$\tau \ \leq (m+n)^{1/\delta}$ in (1.1.3)
\subhead (8.1) Proof of Theorem 1.3 under the additional assumption that $\tau$ 
is bounded by a polynomial in $m+n$ \endsubhead 
All constants implicit in the ``$O$'' and ``$\Omega$'' notation below may depend only on the 
parameter $\delta$. We say that $m$ and $n$ are {\it sufficiently large} provided 
$m+n \geq \gamma(\delta)$ for some constant $\gamma(\delta)>0$. 

As in Corollary 2.2, we represent 
the number $\#(R, C)$ of tables as the integral 
$$\#(R, C)={e^{g(Z)} \over (2 \pi)^{m+n-1}} \int_{\Pi_0} F(s, t) \ ds dt.$$
Let $\UU \subset \Pi_0$ be the neighborhood of the origin as defined in Theorems 6.1 and 7.1.
From Parts (2), (3) and (4) of Theorem 6.1 we conclude that the integrals of $F(s, t)$ and $|F(s, t)|$
over $\UU$ are 
of the same order, that is 
$$\int_{\UU} |F(s, t)| \ ds dt \ \leq \ O\left( \left| \int_{\UU} F(s, t) \ ds dt \right|\right)$$
provided $m+n$ is sufficiently large.
Theorem 7.1 implies then that the integral of $F(s, t)$ over $\Pi_0 \setminus \UU$ is asymptotically 
negligible: for any $\kappa >0$ we have 
$$\left| \int_{\Pi_0 \setminus \UU} F(s, t) \ ds dt \right| \ \leq \ (m+n)^{-\kappa} 
\left| \int _{\UU} F(s, t) \ ds dt \right|$$
provided $m+n > \gamma(\delta, \kappa)$ for some $\gamma(\delta, \kappa) >0$.

We use Part (3) of Theorem 6.1 to compute 
$$\int_{\UU} F(s, t) \ ds dt.$$
Identifying ${\Bbb R}^{m+n-1}$ with the hyperplane $\tau_n=0$ in ${\Bbb R}^{m+n}$, we note that 
$$\Xi=\int_{{\Bbb R}^{m+n-1}} e^{-q} \ ds dt ={\pi^{(m+n-1)/2} \over \sqrt{\det q| {\Bbb R}^{m+n-1}}},$$
and that by Lemma 3.5 we have 
$$\det q| {\Bbb R}^{m+n-1} ={1 \over m+n} \det q|H,$$
where $H$ is the hyperplane orthogonal to the null-space of $q$. 

To conclude the proof, we note that by Lemma 3.1 the values 
of 
$$\mu=\EE f^2 \quad \text{and} \quad \nu=\EE h$$
can be computed with respect to the Gaussian probability measure with the density proportional 
to $e^{-q}$ in an arbitrary hyperplane $L \subset {\Bbb R}^{m+n}$ not containing the 
null-space of $q$.
{\hfill \hfill \hfill} \qed
\bigskip
To handle the case of super-polynomial $\tau$, we use a result of \cite{D+97, Lemma 3}, 
which shows that $\#(R, C) \approx \vl P(R, C)$ provided the margins $R$ and $C$ are large 
enough (it suffices to have $\tau \geq (mn)^2$). Then we note that 
$$\vl P(\alpha R, \alpha C) =\alpha^{(m-1)(n-1)} \vl P(R, C) \quad \text{for} \quad \alpha >0$$
and show that the formula of Theorem 1.3 scales similarly. In the next three lemmas we 
show that the typical matrix of $(\alpha R, \alpha C)$ is approximately the typical matrix 
of $(R, C)$ multiplied by $\alpha$ and that the typical matrix of $(R, C)$ is approximately 
the typical matrix of $(R', C')$ if $R' \approx R$ and $C' \approx C$. 
We then complete our proof of Theorem 1.3.

In Lemmas 8.2 and 8.3 below, all implicit constants in the ``$O$'' notation are absolute.  

\proclaim{(8.2) Lemma} Let $R=\left(r_1, \ldots, r_m\right)$ and $C=\left(c_1, \ldots, c_n\right)$
be positive (not necessarily integer) vectors such that  
$r_1 + \ldots + r_m = c_1 + \ldots + c_n$
and let $Z=\left(\zeta_{jk}\right)$ be the typical matrix maximizing the the value of 
$$g(X)=\sum \Sb 1 \leq j \leq m \\ 1 \leq k \leq n \endSb  \Bigl( \left(x_{jk}+1\right) \ln 
\left(x_{jk}+1 \right) - x_{jk} \ln x_{jk} \Bigr)$$ on the polytope $P(R, C)$ of $m \times n$ non-negative
matrices with row sums $R$ and column sums $C$.

Let 
$$\split &r_-=\min_{j=1, \ldots, m} r_j, \quad c_-=\min_{k=1, \ldots, n} c_k \quad \text{and} \\
&r_+ =\max_{j=1, \ldots, m} r_j, \quad c_+ =\max_{k=1, \ldots, n} c_k. \endsplit$$
Then 
$$\zeta_{jk} \ \geq \ {r_- c_- \over r_+ m} \quad \text{and} \quad \zeta_{jk} \ \geq \ {c_- r_- \over c_+ n} 
\quad \text{for all} \quad j, k.$$
\endproclaim 
\demo{Proof} This is Part (1) of Theorem 3.5 (Theorem 3.3 of the journal version) of \cite{B+08}. 
{\hfill \hfill \hfill} \qed 
\enddemo

\proclaim{(8.3) Lemma} Let $Z=\left(\zeta_{jk}\right)$ be the $m \times n$ typical matrix of margins 
$(R, C)$ such that 
$$\delta \tau \ \leq \ \zeta_{jk} \ \leq \tau \quad \text{for all} \quad j, k,$$
for some $0 < \delta < 1/2$ and some $\tau>1$. 

Let $0 < \alpha < 1$ and let 
$X=\left(\xi_{jk}\right)$ be the typical matrix of margins $(\alpha R, \alpha C)$.

Then the following holds:
\roster
\item We have
$$\xi_{jk} \ \geq \ \alpha \delta^2 \tau \quad \text{for all} \quad j, k.$$
\item Suppose $\alpha \delta^2 \tau > 1$. Then
$$ \Bigl| g(Z) + mn \ln \alpha - g(X)  \Bigr| =O\left({mn \over \alpha \delta^2 \tau}\right).$$
\item There exists an absolute constant $\gamma > 1$ such that if 
$\alpha \delta^4 \tau \geq \gamma m n$ then 
$$\Bigl| \xi_{jk} - \alpha \zeta_{jk} \Bigr| \ \leq \ \epsilon \alpha \zeta_{jk}
 \quad \text{for all} \quad j, k$$
and 
$$ \epsilon =O\left( \sqrt{mn \over \alpha \delta^4 \tau}\right).$$
\endroster
\endproclaim 
\demo{Proof} 
Let $R=\left(r_1, \ldots, r_m\right)$ and $C=\left(c_1, \ldots, c_n \right)$. 
Thus in Lemma 8.2 we have 
$$r_- \ \geq \ \delta n \tau, \quad c_- \ \geq \ \delta m \tau \quad \text{and} \quad 
r_+ \ \leq \ \tau n. $$
Applying Lemma 8.2 to the scaled margins $(\alpha R, \alpha C)$, we obtain Part (1).

Since $\alpha Z \in P(\alpha R, \alpha C)$ and $\alpha^{-1} X \in P(R, C)$, we have 
$$g(Z) \ \geq \ g\left(\alpha^{-1} X \right) \quad \text{and} \quad g(X) \ \geq \ g(\alpha Z). \tag8.3.1$$
Since for $x \geq 1$ we have 
$$\aligned g(x)=&(x+1)\ln (x+1) - x\ln x +(x+1) \ln x - (x+1) \ln x \\=&(x+1) \ln \left(1+{1 \over x} \right) +\ln x
 = 1+ \ln x + O \left({1 \over x} \right), \endaligned \tag8.3.2$$ 
from (8.3.1) we obtain Part (2).

Let us consider the interval $[X, \alpha Z] \subset P(\alpha R, \alpha C)$. Since 
$g$ is concave, we have 
$$\bigl| g(Y) - g(X) \bigr|  =O\left( {mn \over \alpha \delta^2 \tau}\right)
\quad \text{for all} \quad Y \in [X, \alpha Z].\tag8.3.3$$

Suppose that for some $0 < \epsilon < 1/2$ we have 
$$\left| \xi_{jk} - \alpha \zeta_{jk}\right| > \epsilon \alpha \zeta_{jk} \quad \text{for some} \quad j, k.$$
Then there is a matrix $Y \in [X, \alpha Z]$, $Y=\left(\eta_{jk}\right)$, such that 
$\left| \eta_{jk} - \alpha \zeta_{jk} \right| = \epsilon \alpha \zeta_{jk}$. 
We note that 
$$g''(x) =-{1 \over x(x+1)}$$
and, in particular, 
$$g''(x) \ \leq \ -{1 \over 8 \alpha^2  \tau^2} \quad \text{for all} \quad x \in
 \left[ \eta_{jk},\ \alpha \zeta_{jk} \right].$$
Next, we are going to exploit the strong concavity of $g$ and use the following standard inequality:
if $g''(x) \leq -\beta$ for some $\beta >0$ and all $ a \leq x \leq b$ then 
$$g\left({a+ b \over 2} \right) -{1 \over 2} g(a)-{1 \over 2} g(b) \ \geq \ {\beta(b-a)^2 \over 8}.$$
Applying the above inequality to $g$ with $a=\eta_{jk}$, $b=\alpha \zeta_{jk}$ and 
$\beta=1/8 \alpha^2 \tau^2$, we obtain 
$$g\left({\eta_{jk} + \alpha \zeta_{jk} \over 2} \right) -
{1 \over 2} g\left(\eta_{jk}\right) - {1 \over 2} g \left(\alpha \zeta_{jk} \right) \ \geq \ 
{\epsilon^2  \delta^2   \over 64}. $$
Let $W=(Y+\alpha Z)/2$. Then $W \in [X, \alpha Z]$ and by (8.3.3)
$$g(W) \ \geq \ g(X) + {\epsilon^2 \delta^2 \over 64} - O\left({mn \over \alpha \delta^2 \tau} \right).$$
Since $g(W) \leq g(X)$, the proof follows.
{\hfill \hfill \hfill} \qed
\enddemo

\proclaim{(8.4) Lemma} Let $Z=\left(\zeta_{jk}\right)$ be the $m \times n$ typical matrix of margins 
$R=\left(r_1, \ldots, r_m\right)$ and $C=\left(c_1, \ldots, c_n \right)$ such that 
$$\delta \tau \ \leq \ \zeta_{jk} \ \leq \ \tau \quad \text{for all} \quad j, k,$$
for some $0 < \delta < 1/2$ and some $\tau >1$. 

Let $0 < \epsilon < 1/2$ and let 
$X=\left(\xi_{jk} \right)$ be the typical matrix of some margins $R=\left(r_1', \ldots, r_m' \right)$
and $C'=\left(c_1', \ldots, c_n' \right)$ such that 
$$\split &(1-\epsilon) r_j \ \leq \ r_j' \ \leq \ r_j \quad \text{for} \quad j=1, \ldots, m \qquad 
\text{and} \\
&(1-\epsilon) c_k \ \leq \ c_k' \ \leq \ c_k \quad \text{for} \quad k=1, \ldots, n. \endsplit$$
Suppose that $\delta^2 \tau >1$. 
Then the following holds:
\roster
\item We have
$$\Big| g(X) - g(Z) \Big| = O\left(mn \epsilon \right).$$
\item There exists an absolute constant $\gamma$ such that if $\epsilon \leq \gamma \delta^2/mn$
then 
$$\big| \xi_{jk} - \zeta_{jk} \big| \ \leq \ \beta \zeta_{jk} \quad \text{for all} \quad j, k$$
and 
$$\beta=O\left(\sqrt{mn \epsilon \over \delta^2}\right).$$
\endroster
\endproclaim 
\demo{Proof} Let $A=\left(a_1, \ldots, a_m\right)$ and $B=\left(b_1, \ldots, b_n\right)$ be margins and let 
$A'=\left(a_1', \ldots, a_m' \right)$ and $B'=\left(b_1', \ldots, b_n'\right)$ be some other margins such that 
$a_j' \leq a_j$ and $b_k' \leq b_k$ for all $j$ and $k$. Then there exists a non-negative 
matrix $D$ with margins $a_j-a_j'$ and $b_k-b_k'$ and for such a $D$ we have 
$Y + D \subset P(A, B)$ for all $Y \in P(A', C')$. 
Snce $g$ is monotone increasing, we obtain 
$$\split g\bigl((1-\epsilon) Z \bigr) \ &\leq \ \max_{Y \in (1-\epsilon) P(R, C)} g(Y) \ \leq \ 
g(X) =\max_{Y \in P(R', C')} g(Y) \\  &\leq \ \max_{Y \in P(R, C)} g(Y) =g(Z). \endsplit$$
Hence 
$$g\bigl((1-\epsilon) Z \bigr) \ \leq \ g(X) \ \leq \ g(Z)$$
and from (8.3.2) we deduce Part (2). 

We note that $Z$ is the maximum point of $g$ on the polytope of non-negative $m \times n$ matrices 
with the row 
sums not exceeding $R$ and column sums not exceeding $C$. Therefore, 
$$\Big| g(Y)- g(Z) \Big| = O(mn \epsilon) \quad \text{for all} \quad Y \in [X, Z]. \tag8.4.1$$
Suppose that for some $0 < \beta < 1/2$ we have 
$$\big| \xi_{jk} - \zeta_{jk} \big|  > \beta  \zeta_{jk} \quad \text{for some} \quad j, k.$$
Then there is a matrix $Y \in [X, Z]$, $Y=\left(\eta_{jk}\right)$, such that 
$\big| \eta_{jk}- \zeta_{jk} \big| =\beta \zeta_{jk}$. As in the proof of Lemma 8.3, we 
argue that 
$$g''(x) \ \leq \ -{1 \over 8 \tau^2} \quad \text{for all} \quad x \in [\eta_{jk},\ \zeta_{jk}]$$
and that 
$$g\left( {\eta_{jk} + \zeta_{jk} \over 2}\right)-{1 \over 2} g\left(\eta_{jk}\right) -{1 \over 2} g\left(\zeta_{jk}\right) 
\ \geq \ {\beta^2 \delta^2 \over 64}.$$
Let $W=(Y + Z)/2$. Then $W \in [Y, Z]$ and by (8.4.1) 
$$g(W) \ \geq \ g(Z) +{\beta^2 \delta^2 \over 64} - O(mn \epsilon).$$
Since $g(W) \leq g(Z)$, the proof follows.
{\hfill \hfill \hfill} \qed
\enddemo

\subhead (8.5) Proof of Theorem 1.3 \endsubhead 
All implicit constants in the ``$O$'' and ``$\Omega$'' notation below may depend on parameter $\delta$ only.

In view of Section 8.1, without loss of generality we assume that 
$\tau \ \geq \ (m+n)^{10}$ in (1.1.3). 
As follows by \cite{D+97}, as long as $\tau \geq (m+n)^2$ we have 
$$\#(R, C) =\vl P(R, C)\left(1 + O\left({1 \over m+n} \right)\right),$$
where $\vl P(R, C)$ is the volume of the polytope of the set of $m \times n$ non-negative
matrices with row sums $R$ and column sums $C$ normalized in such a way that the 
volume of the fundamental domain of the $(m-1)(n-1)$-dimensional lattice consisting 
of the $m \times n$ integer matrices with zero row and column sums is equal to 1.

Let $\alpha =(m+n)^{9} \tau^{-1}$ 
and let 
$$\hat{R}=\left(\hat{r}_1, \ldots, \hat{r}_m \right) \quad \text{and} \quad 
\hat{C}=\left(\hat{c}_1, \ldots, \hat{c}_n \right)$$ 
be positive integer margins (so $\hat{r}_1 + \ldots + \hat{r}_m =\hat{c}_1 + \ldots + \hat{c}_n$) such that 
$$\split &(1-\epsilon) \alpha r_j \ \leq \ \hat{r}_j \ \leq \ \alpha r_j \quad \text{and} \quad (1-\epsilon) \alpha c_k \ \leq \hat{c}_k \ \leq \ 
\alpha c_k \\ &\qquad \qquad  \text{for some} \quad 0 < \epsilon < (m+n)^{-7}. \endsplit$$
Then 
$$ \aligned  \#(R, C) \approx &\vl P(R, C) 
\approx \alpha^{(m-1)(1-n)} \vl P(\hat{R}, \hat{C}) \\ \approx &\alpha^{(m-1)(1-n)} \#(\hat{R}, \hat{C}), \endaligned
 \tag8.5.1$$ 
where ``$\approx$'' denotes the equality up to a $O\left((m+n)^{-1}\right)$ relative error.

Let $Z=\left(\zeta_{jk}\right)$ be the typical matrix of margins $(R, C)$ and let 
$\hat{Z}=\left(\hat{\zeta}_{jk}\right)$ be the typical matrix of margins $(\hat{R}, \hat{C})$. 
By Lemmas 8.3 and 8.4, we have 
$$\aligned &\big| g(Z) + mn \ln \alpha - g(\hat{Z}) \big| =O \left({1 \over (m+n)^{5}} \right) \quad \text{and} \\
&\big| \hat{\zeta}_{jk} - \alpha \zeta_{jk} \big| = O \left({ \alpha \zeta_{jk} \over (m+n)^2}\right) 
\quad \text{for all} \quad j, k. \endaligned \tag8.5.2$$

Let $q, \hat{q}: {\Bbb R}^{m+n} \longrightarrow {\Bbb R}$ be the quadratic forms associated by (1.2.1)
with margins $(R, C)$ and $(\hat{R}, \hat{C})$ respectively. Then by the second estimate 
of (8.5.2) it follows that 
$$\hat{q}(s, t) \approx \alpha^2 q(s, t) \quad \text{for all} \quad (s, t) \in {\Bbb R}^{m+n},$$ 
where ``$\approx$'' stands for the equality up to a $O\left((m+n)^{-2}\right)$ relative error. 
It follows then by the first estimate of (8.5.2) that the Gaussian term (1.3.1) for margins $(R, C)$, 
up to a relative error of $O\left( (m+n)^{-1}\right)$, is obtained by multiplying the 
Gaussian term for margins 
$(\hat{R}, \hat{C})$ by $\alpha^{(m-1)(1-n)}$. 

Similarly, we show that the 
Edgeworth correction factor (1.3.2) changes negligibly as we pass from $(R, C)$ to $(\hat{R}, \hat{C})$.
By making substitutions 
$$(s, t) \longmapsto \tau^{-1} (s, t) \quad \text{and} \quad (s, t) \longmapsto \alpha^{-1} \tau^{-1} (s, t)$$
respectively, we express the quantities $(\mu, \nu)$ for margins $(R, C)$ and 
$(\hat{\mu}, \hat{\nu})$ for margins 
$(\hat{R}, \hat{C})$ as 
$$\mu = \EE f^2, \ \nu = \EE h \quad \text{and} \quad \hat{\mu} =\EE \hat{f}^2, \ \hat{\nu} =\EE \hat{h},$$
where the expectations $\mu$ and $\nu$ are taken with respect to the Gaussian measure 
on $H$ with the density proportional to $e^{-\psi}$ and the expectations $\hat{\mu}$ and 
$\hat{\nu}$ are taken with respect to the Gaussian measure with the density proportional to 
$e^{-\hat{\psi}}$, where $\psi$ and $\hat{\psi}$ are positive definite quadratic forms within 
a relative error of $O\left((m+n)^{-2}\right)$ of each other. Moreover, $f^2$ and $\hat{f}^2$ 
are homogeneous polynomials of degree 6 and $h$ and $\hat{h}$ are homogeneous polynomials 
of degree 4 such that 
$$\split &f^2(s, t), \  \hat{f}^2(s, t) =
O\left( \sum_{j, k} \left| s_j + t_k\right|^3 \right)^2, \\
&h(s, t), \ \hat{h}(s, t) =O\left( \sum_{j, k} \left(s_j + t_k \right)^4 \right) \qquad \text{and} \\
 &\left| f^2(s, t) - \hat{f}^2(s, t)\right| = O\left({1\over (m+n)^2}\right) 
\left(\sum_{j, k} \left|s_j + t_k \right|^3 \right)^2, \\
&\left| h(s, t) - \hat{h}(s, t) \right| =O\left({1\over (m+n)^2}\right) 
\sum_{j, k}   \left(s_j + t_k \right)^4. \endsplit$$
Since by Lemma 3.6, the minimum eigenvalues of $\psi$ and $\hat{\psi}$ are $\Omega(m+n)$,
standard estimates imply that $\exp\left\{-\mu/2 + \nu \right\}$ approximates 
$\exp\left\{-\hat{\mu}/2 +\hat{\nu}\right\}$ within a $O\left((m+n)^{-1}\right)$ relative error.

We have 
$$\hat{\zeta}_{jk} = O\left((m+n)^{9}\right) \quad \text{for all} \quad j, k,$$
and hence by the result of Section 8.1 we can apply Theorem 1.3 to estimate 
$\#(\hat{R}, \hat{C})$.  The proof then follows by 
(8.5.1).
{\hfill \hfill \hfill} \qed

\enddocument

\end